\documentclass{article}

\pdfoutput=1

\usepackage[T1]{fontenc}
\usepackage[utf8]{inputenc}

\usepackage{hyperref}
    \hypersetup{
        pdftex,
        pdfauthor={Gustav Ludvigsson, Kyle R. Steffen, Simon Sticko, Siyang Wang, Qing Xia, Yekaterina Epshteyn, Gunilla Kreiss},
        pdftitle={High-order numerical methods for 2D parabolic problems in single and composite domains},
        pdfsubject={Original research},
        pdfkeywords={parabolic problems, interface models, level set, complex geometry, discontinuous solutions, SBP--SAT finite difference, difference potentials, spectral approach, finite element method, cut elements, immersed boundary, stabilization, higher order accuracy and convergence},
        pdfproducer={Latex with hyperref},
        pdfcreator={pdflatex},
        colorlinks,
        linkcolor=blue,
        citecolor=blue,
        urlcolor=blue,
        final
    } 
\usepackage{microtype}
\usepackage{doi}
\usepackage{booktabs,siunitx} 

\usepackage{authblk}

\usepackage{afterpage}
\usepackage{subcaption,float} 
    \floatstyle{plaintop} \restylefloat{table} 

\usepackage{hyphenat}

\usepackage{graphicx}

\usepackage{amsmath,amssymb,amsthm,mathtools}

\DeclarePairedDelimiter{\abs}{\lvert}{\rvert}
\newcommand{\RR}{\mathbb{R}}
\newcommand{\Tr}{\operatorname{Tr}}

\newtheorem{conj}{Conjecture}
\newtheorem{definition}{Definition}
\newtheorem{remark}{Remark}
\newtheorem{theorem}{Theorem}

\newcommand{\siyang}[3]{#2}
\newcommand{\redtextst}[3]{#2}
\newcommand{\redtext}[2]{#1}

\title{High-order numerical methods for 2D parabolic problems in single and composite domains
}

\author[1]{Gustav Ludvigsson\thanks{gustav.ludvigsson@it.uu.se}}
\author[2]{Kyle R. Steffen\thanks{steffen@math.utah.edu}}
\author[1]{Simon Sticko\thanks{simon.sticko@it.uu.se}}
\author[1]{Siyang Wang\thanks{siyang.wang@it.uu.se}}
\author[2]{Qing Xia\thanks{xia@math.utah.edu}}
\author[2]{Yekaterina Epshteyn\thanks{epshteyn@math.utah.edu}}
\author[1]{Gunilla Kreiss\thanks{gunilla.kreiss@it.uu.se}}
\affil[1]{\small Department of Information Technology, Uppsala University, Box 337, 751 05 Uppsala, Sweden}
\affil[2]{\small Department of Mathematics, The University of Utah, 155 S 1400 E Rm. 233, Salt Lake City, UT 84112, USA}

\begin{document}


\maketitle

\begin{abstract}
In this work, we discuss and compare three methods for the numerical approximation of constant- and variable-coefficient diffusion equations in both single and composite domains with possible discontinuity in the solution/flux at interfaces, considering (i) the Cut Finite Element Method; (ii) the Difference Potentials Method; and (iii) the summation--by--parts Finite Difference Method. First we give a brief introduction for each of the three methods. Next, we propose benchmark problems, and consider numerical tests--with respect to accuracy and convergence--for linear parabolic problems on a single domain, and continue with similar tests for linear parabolic problems on a composite domain (with the interface defined either explicitly or implicitly). Lastly, a comparative discussion of the methods and numerical results will be given. 


\end{abstract}



\section{Introduction} \label{sec:Introduction}


Designing methods for the high-order accurate numerical approximation of partial differential equations (PDE) posed on composite domains with interfaces, or on irregular and geometrically complex domains, is crucial in the modeling and analysis of problems from science and engineering. Such problems may arise, for example, in materials science (models for the evolution of grain boundaries in polycrystalline materials), fluid dynamics (the simulation of homogeneous or multi-phase fluids), engineering (wave propagation in an irregular medium or a composite medium with different material properties), biology (models of blood flow or the cardiac action potential), etc. The analytic solutions of the underlying PDE may have non-smooth or even discontinuous features, particularly at material interfaces or at interfaces within a composite medium. Standard numerical techniques involving finite-difference approximations, finite-element approximation, etc., may fail to produce an accurate approximation near the interface, leading one to consider and develop new techniques.


There is extensive existing work addressing numerical approximation of PDE posed on composite domains with interfaces or irregular domains, for example, the 
boundary integral method \cite{bouchon1989boundary,mayo1984fast}, 
difference potentials method \cite{phdthesisJA,albright2015high,Epshteyn2014,epshteyn2015solution,MedvinskyTsynkovTurkel2012,Ryabenkii2002}, 
immersed boundary method \cite{fadlun2000combined,kim2001immersed,peskin2002immersed,tseng2003ghost}, 
immersed interface method \cite{adams2002immersed,leveque1994immersed,leveque1997immersed,linnick2005high,sethian2000structural}, 
ghost fluid method \cite{fedkiw1999non,gibou2005fourth,liu2003ghost,liu2003convergence},  
the matched interface and boundary method \cite{xia2011matched,yu2007three,yu2007matched,zhou2006high}, 
Cartesian grid embedded boundary method \cite{crockett2011cartesian,johansen1998cartesian,mccorquodale2001cartesian,ye1999accurate}, 
multigrid method for elliptic problems with discontinuous coefficients on an arbitrary interface \cite{coco2012second}, 
virtual node method \cite{bedrossian2010second,hellrung2012second}, 
Voronoi interface method \cite{guittet2015solving,guittet2017voronoi}, 
\siyang{summation--by--parts}{the finite difference method \cite{appelo2009stable,Berg2012,virta2013high,Virta2014,wang2016high,wang2017fourth}  and finite volume method \cite{Demirdzic1994,Gong2013} based on mapped grids,}{Comment 11 of reviewer 1} 
or cut finite element method \cite{cutfem_2015,burman_hansbo_cut_II,cutfem_coupled_bulk_surface_2016,hansbo_unfitted_2002,hansbo_cut_2014,stickoLowerOrder2016,wadbro2013uniformly}. 
\redtextst{In spite of}{Indeed, there have been}{*} great advances in numerical methods for the approximation of PDE posed on composite domains with interfaces, or on irregular domains. \redtext{However,}{*} it is still a challenge to design high-order accurate and \redtextst{efficient}{computationally-efficient}{*} methods \redtext{for PDE posed in these complicated geometries}{*}, especially for time-dependent problems, \redtext{problems with variable coefficients, or problems with general boundary/interface conditions.}{* $=$ Comment 2 of Reviewer 1}

The aim of this work is to establish benchmark (test) problems for the numerical approximation of parabolic PDE defined in irregular or composite domains. 
\redtext{The considered models (Section~\ref{sec:Intro:ContinuousProblem}) arise in the study of mass or heat diffusion in single or composite materials, or as simplified models in other areas (\textit{e.g.}, biology, materials science, etc.).}{Comment 3 Reviewer 1} 
\redtextst{In particular, the formulated}{The formulated test}{Comment 3 Reviewer 1} problems (Section~\ref{sect:TestProblems}) are intended (a) to be suitable for comparison of high-order accurate numerical methods -- and will be used as such in this study -- and (b) to be useful in further research. 
Moreover, the proposed problems include a wide variety of possibilities relevant in applications, which any robust numerical method should resolve accurately, including constant diffusion; time-varying diffusion; high frequency oscillations in the analytical solution; large jumps in diffusion coefficients, solution, and/or flux; etc. 
For now, we will consider a simplified geometrical setting, with the intent of setting a ``baseline'' from which further research, or more involved comparisons, might be conducted.  Therefore, in Section~\ref{sec:Intro:ContinuousProblem} we will introduce two circular geometries, which are defined either explicitly, or implicitly via a level set function. 

In Section~\ref{sec:Intro:OverviewMethods}, we briefly introduce the numerical methods we will consider in this work, \textit{i.e.}, second- and fourth-order versions of (i) the Cut Finite Element Method (cut--FEM); (ii) the Difference Potentials Method (DPM), with Finite Difference approximation as the underlying discretization in the current work; and (iii) the summation--by--parts Finite Difference Method combined with the simultaneous approximation term technique (SBP--SAT--FD). 
These three methods are all modern numerical methods which may be designed for problems in irregular or composite domains, allowing for high-order accurate numerical approximation, even at points close to irregular interfaces or boundaries. 
We will apply each method to the formulated benchmark problems, and compare results. From the comparisons, we expect to learn what further developments of the methods at hand would be most important. 

To resolve geometrical features of irregular domains, both cut--FEM and DPM use a Cartesian grid on top of the domain, which need not conform with boundaries or interfaces.  These types of methods are often characterized as ``immersed'' or ``embedded''.  In the finite difference framework, embedded methods for parabolic problems are developed in \cite{Abarbanel1997,Ditkowski2009}.   \siyang{For comparison with cut--FEM and DPM, however, in this paper we use a conforming approach based on the finite difference method -- the SBP--SAT--FD method, which resolves geometrical features by curvilinear mapping.}{For comparison with cut--FEM and DPM, however, in this paper we use a finite difference method based on a conforming approach. The finite difference operators we use satisfy a summation--by--parts principle. Then, in combination with the SAT method to weakly impose boundary and interface conditions, an energy estimate of the semi--discretization can be derived to ensure stability. In addition, we use curvilinear grids and transfinite interpolation to resolve complex geometries.}{Comment 10 of reviewer 1}

For recent work on SBP--SAT--FD for wave equations in composite domains, see \cite{Berg2012,Carpenter2009,Virta2014,wang2016high}, and the two review papers \cite{Del2014Review,Svard2014}; for recent work in DPM for elliptic/parabolic problems in composite domains with interface defined explicitly, see \cite{phdthesisJA,albright2017high,AlbrightEpshteynSteffen2015,albright2015high,Epshteyn2014,epshteyn2015solution,epshteyn2015high,MedvinskyTsynkovTurkel2012,medvinsky2016solving,Ryabenkii2002}; and for recent work in cut--FEM see \cite{cutfem_2015,burman_hansbo_cut_II,cutfem_coupled_bulk_surface_2016,hansbo_unfitted_2002,massing_stabilized_2014,stickoLowerOrder2016}. 

The paper is outlined as follows. In Section~\ref{sec:Intro:ContinuousProblem}, we give brief overview of the continuous formulation of the parabolic problems in a single domain or a composite domain.  In Section~\ref{sec:Intro:OverviewMethods}, we give introductions to the basics of the three proposed methods: cut--FEM, DPM, and SBP--SAT--FD. In Section~\ref{sect:TestProblems}, we formulate the numerical test problems. In Section~\ref{sec:NumericalResults}, we present extensive numerical comparisons of errors and convergence rates, between the second- and fourth-order versions of each method. The comparisons include single domain problems with constant or time-dependent diffusivity; and interface problems with interface defined explicitly, or implicitly by a level set function. In Section~{\ref{sec:Discussion}}, we give a comparative discussion of the three methods and the numerical results, together with a discussion on future research directions. Lastly, in Section~\ref{sec:Conclusions}, we give our concluding remarks. 

\section{Statement of problem} \label{sec:Intro:ContinuousProblem}

\afterpage{
\begin{figure}[tb]
    \centering
    \begin{subfigure}[b]{0.5\linewidth}
        \centering
        \includegraphics[width=2in]{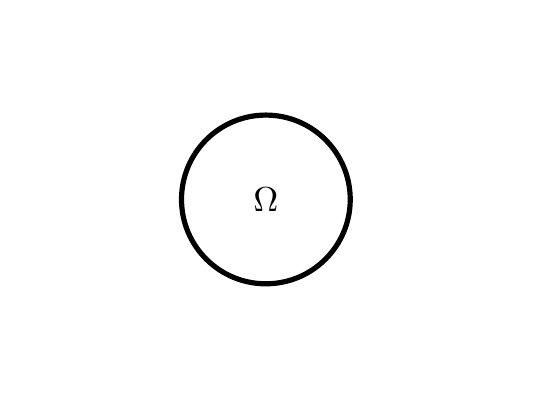}
        \caption{\label{fig:SingDomainSubFig}}
    \end{subfigure}%
    \begin{subfigure}[b]{0.5\linewidth}
        \centering
        \includegraphics[width=2in]{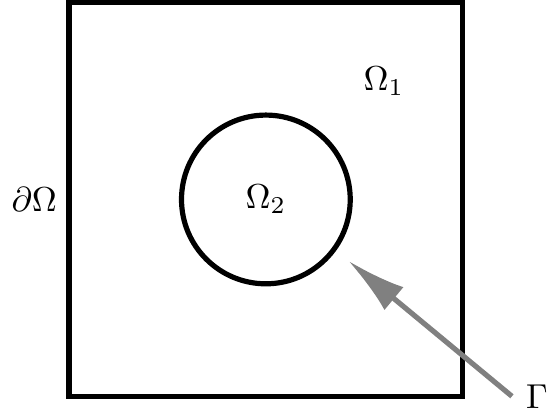}
        \caption{\label{fig:CompDomainSubFig}}
    \end{subfigure}
    \caption{The \textbf{(a)} single domain $\Omega$ and \textbf{(b)} composite domain $\Omega = \Omega_1 \cup \Omega_2$.  In \textbf{(b)}, $\partial\Omega_1$ has two connected components: the boundary $\partial\Omega$ and interface $\Gamma = \partial\Omega_2$.}
    \label{fig:CompDomain}
\end{figure}%
}%
In this section, we describe two diffusion problems, which will be the setting for our proposed benchmark (test) problems in Section~\ref{sect:TestProblems}.  \redtext{(Recall from Section~\ref{sec:Introduction} that these models arise, for example, in the study of mass or heat diffusion.)}{Comment 3 Reviewer 1}  For brevity, in the following discussion, we denote $u := u(x, y, t)$ and $u_s := u_s(x, y, t)$, with $s = 1, 2$. 
\subsection{The single domain problem}
  First, we consider the linear parabolic PDE on a single domain $\Omega$ (\textit{e.g.}, Figure~\ref{fig:SingDomainSubFig}), with variable diffusion $\lambda(t)$: %
  \begin{align}
    \dfrac{\partial u}{\partial t} &= \nabla \cdot (\lambda(t) \nabla u) + f(x, y, t), \ \  (x, y, t) \in \Omega \times (0, T], \label{eqn:PDESD1} \\
    \intertext{subject to initial and Dirichlet boundary conditions:}
    u(x, y, 0) &= u^0(x, y), \ \  (x, y) \in \Omega \quad \mbox{and} \quad u = \psi(x, y, t), \ \  (x, y, t) \in \partial\Omega \times (0, T]. \label{eqn:PDESD2-3} 
  \end{align}%
  %
  Here, the initial and boundary data $u^0(x, y)$ and $\psi(x, y, t)$, the diffusion coefficient $\lambda(t)$, the forcing function $f(x,y,t)$, and the final time $T$ are known (given) data.

\subsection{The composite domain problem}
  Next, we consider the linear parabolic PDE on a composite domain $\Omega := \Omega_1 \cup \Omega_2$ (\textit{e.g.}, Figure~\ref{fig:CompDomainSubFig}), with constant diffusion coefficients $(\lambda_1, \lambda_2)$: %
  \begin{align}
    \dfrac{\partial u_1}{\partial t} &= \nabla \cdot (\lambda_1 \nabla u_1) + f_1(x, y, t), \ \  (x, y, t) \in \Omega_1 \times (0, T], \label{eqn:PDECD1} \\
    \dfrac{\partial u_2}{\partial t} &= \nabla \cdot (\lambda_2 \nabla u_2) + f_2(x, y, t), \ \  (x, y, t) \in \Omega_2 \times (0, T], \label{eqn:PDECD2} \\
    \intertext{subject to initial conditions:}
    u_1(x, y, 0) &= u_1^0(x, y), \ \  (x, y) \in \Omega_1, \label{eqn:PDECD3} \\
    u_2(x, y, 0) &= u_2^0(x, y), \ \  (x, y) \in \Omega_2, \label{eqn:PDECD4} \\
    \intertext{Dirichlet boundary conditions:}
    u_1 &= \psi(x, y, t), \ \  (x, y, t) \in \partial\Omega \times (0, T], \label{eqn:PDECD5} \\
    \intertext{and interface/matching conditions:}
    u_1 - u_2 &= \mu_1(x, y, t), \ \  (x, y, t) \in \Gamma \times (0, T], \label{eqn:PDECD6} \\
    \lambda_1 \frac{\partial u_1}{\partial n} - \lambda_2 \frac{\partial u_2}{\partial n} &= \mu_2(x, y, t), \ \  (x, y, t) \in \Gamma \times (0, T].\label{eqn:PDECD7}
  \end{align}%
  In formula \eqref{eqn:PDECD7}, $\frac{\partial u_s}{\partial n}$, $s = 1, 2$ denotes the normal derivative at the interface $\Gamma$, \textit{i.e.}, $\frac{\partial u_s}{\partial n} = \nabla u_s \cdot \mathbf{n}$, where $\mathbf{n}$ is the outward unit normal vector at the interface $\Gamma$.  
  The initial, boundary, and interface data $u_1^0(x, y)$, $u_2^0(x, y)$, $\psi(x, y, t)$, $\mu_1(x, y, t)$, and $\mu_2(x, y, t)$; the diffusion coefficients $(\lambda_1, \lambda_2)$; the forcing functions $f_1(x,y,t)$ and $f_2(x,y,t)$; and the final time $T$ are some known (given) data.
  %


\begin{remark} \label{rmk:MotivationForGeometries}
  We consider the circular geometries depicted in Figure~\ref{fig:CompDomain} as the geometrical setting for our proposed benchmark problems in this work.  In applications (Section~\ref{sec:Introduction}), other geometries will likely be considered, some much more complicated than Figure~\ref{fig:CompDomain}.  While our methods can handle more complicated geometry, this is (to the best of our knowledge) the first work looking to establish benchmarks -- and compare numerical methods -- for parabolic interface problems (\ref{eqn:PDECD1}--\ref{eqn:PDECD7}).  As such, we think that the geometries in Figure~\ref{fig:CompDomain} are a good ``baseline'' -- without all the added complexities that more complicated geometries might produce -- from which further research, or more involved comparisons, might be done.  

  \redtext{To be more specific, we aim to define a simple set of test problems that can be easily implemented and tested for any numerical scheme of interest. With circular domains, it suffices for us to compare/contrast performance of the numerical methods on a simple geometry with smooth boundary \textit{versus} on a composite domain with fixed interface (explicit or implicit). The approximation of the solution to such composite-domain problems are already challenging for any numerical methods, since (i) the solution may fail to be smooth (or may be discontinuous) at the interface, and (ii) there may be discontinuous material coefficients ($\lambda_1 \not= \lambda_2$).}{Comment 1 Reviewer 1}
\end{remark}

\begin{remark}
  For both the single and composite domain problems, we could also consider other boundary conditions, \textit{e.g.}, a Neumann boundary condition as in \cite{albright2015high,cutfem_2015}, etc. 
\end{remark}

  %
  %
  %
  %
  %
  %

\section{Overview of numerical methods} \label{sec:Intro:OverviewMethods}

  \subsection{Cut--FEM} \label{sec:Overview:CutFEM}


\newcommand{\dr}{\mathrm{d}}
\newcommand{\pd}{\partial}
\newcommand{\dd}[2]{\frac{\dr#1}{\dr#2}}
\newcommand{\pdd}[2]{\frac{\partial #1}{\partial #2}}
\newcommand{\pddc}[3]{ \left( \frac{\partial #1}{\partial #2} \right)_#3 }
\newcommand{\ddn}[3]{\frac{\dr^{#1} #2}{\dr #3^{#1}}}
\newcommand{\pddn}[3]{\frac{\partial^{#1} #2}{\partial #3^{#1}}}
\newcommand{\ddop}[1]{\frac{\dr}{\dr #1}}
\newcommand{\pddop}[1]{\frac{\partial}{\partial #1}}
\newcommand{\ddopn}[2]{\frac{\dr^{#1}}{\dr #2^{#1}}}
\newcommand{\pddopn}[2]{\frac{\partial^{#1}}{\partial #2^{#1}}}

\newcommand{\bsprod}[2]{\left\langle #1,#2\right\rangle }
\newcommand{\evaluated}[2]{\left.#1\right|_{#2}}
\newcommand{\T}{\mathcal{T}}
\newcommand{\Fset}{\mathcal{F}}
\newcommand{\M}{\mathcal{M}}

\newcommand{\idx}[2]{{#1}_{#2}}
\newcommand{\pair}[1]{\{ \idx{#1}{1} , \idx{#1}{2} \}}

\newcommand{\ordo}[1]{\mathcal{O}(#1)}
\newcommand{\R}{\mathbb{R}}
\newcommand{\im}{\mathrm{i}}
\newcommand{\e}{\mathrm{e}}

\newcommand{\finalTime}{ t_f }
\newcommand{\elementSize}{h_T}
In this section, we give a brief presentation of the cut--FEM method.  For a more detailed presentation of cut--FEM, see, for example, \cite{cutfem_2015,burman_hansbo_cut_II,massing_stabilized_2014}.

\begin{figure}[tbh] 
\centering
\begin{subfigure}[b]{0.3\linewidth}
\includegraphics[width=\columnwidth]{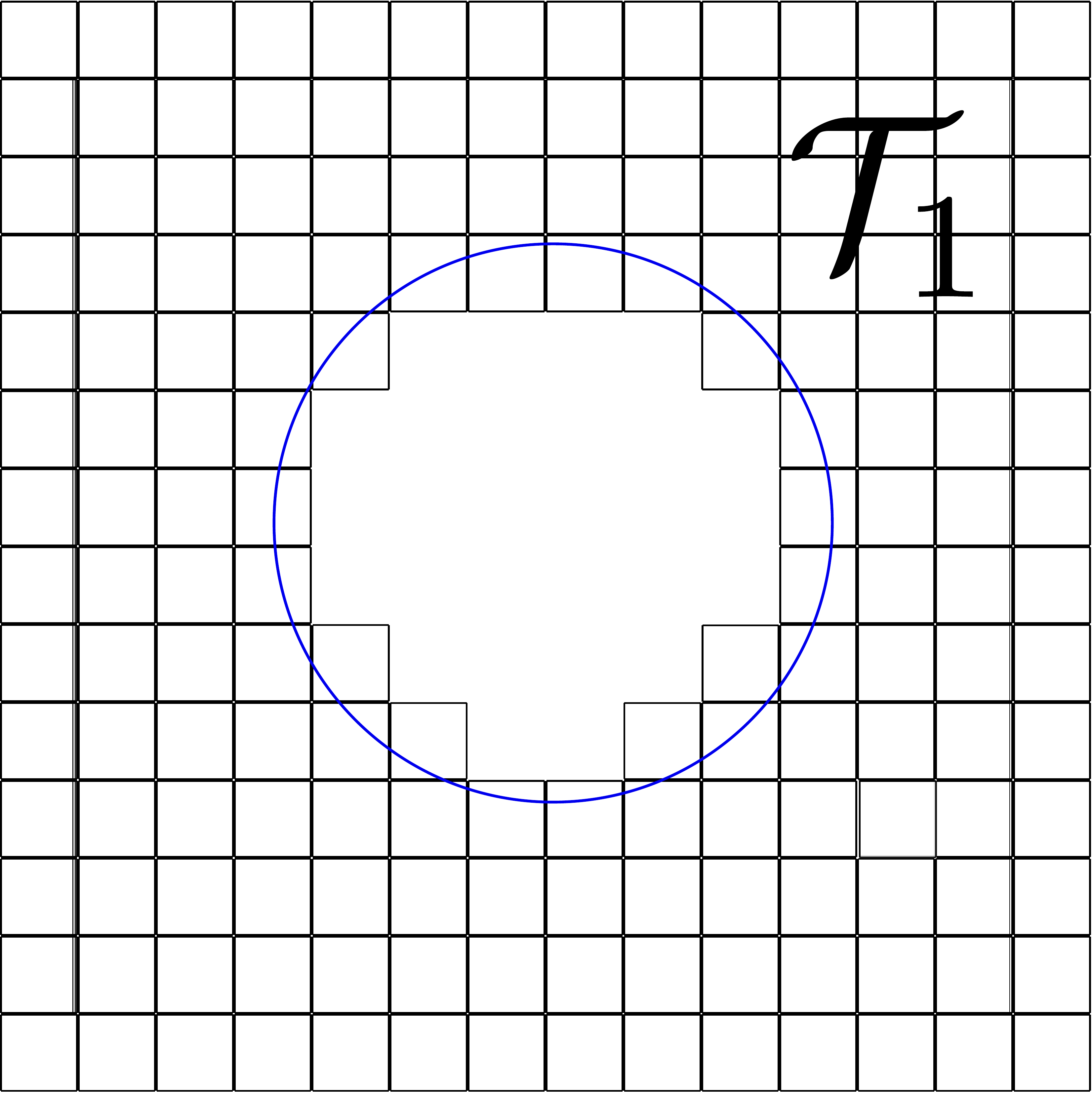}
\caption{\label{fig:outerTriangulation}}
\end{subfigure}
\hfill
\begin{subfigure}[b]{0.33\linewidth}
\includegraphics[width=\columnwidth]{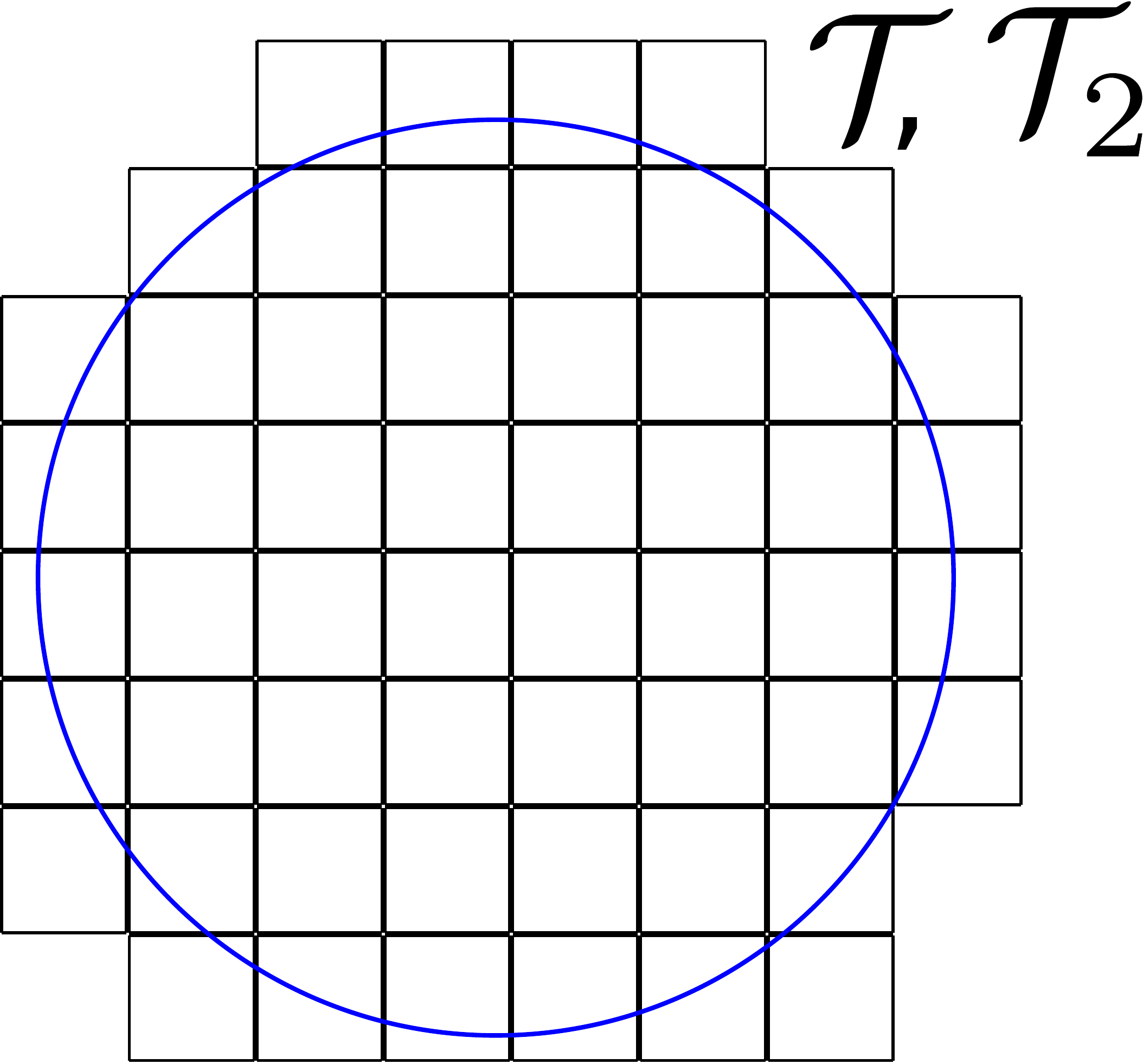}
\caption{\label{fig:innerTriangulation}}
\end{subfigure}
\hfill
\begin{subfigure}[b]{0.3\linewidth}
\includegraphics[width=\columnwidth]{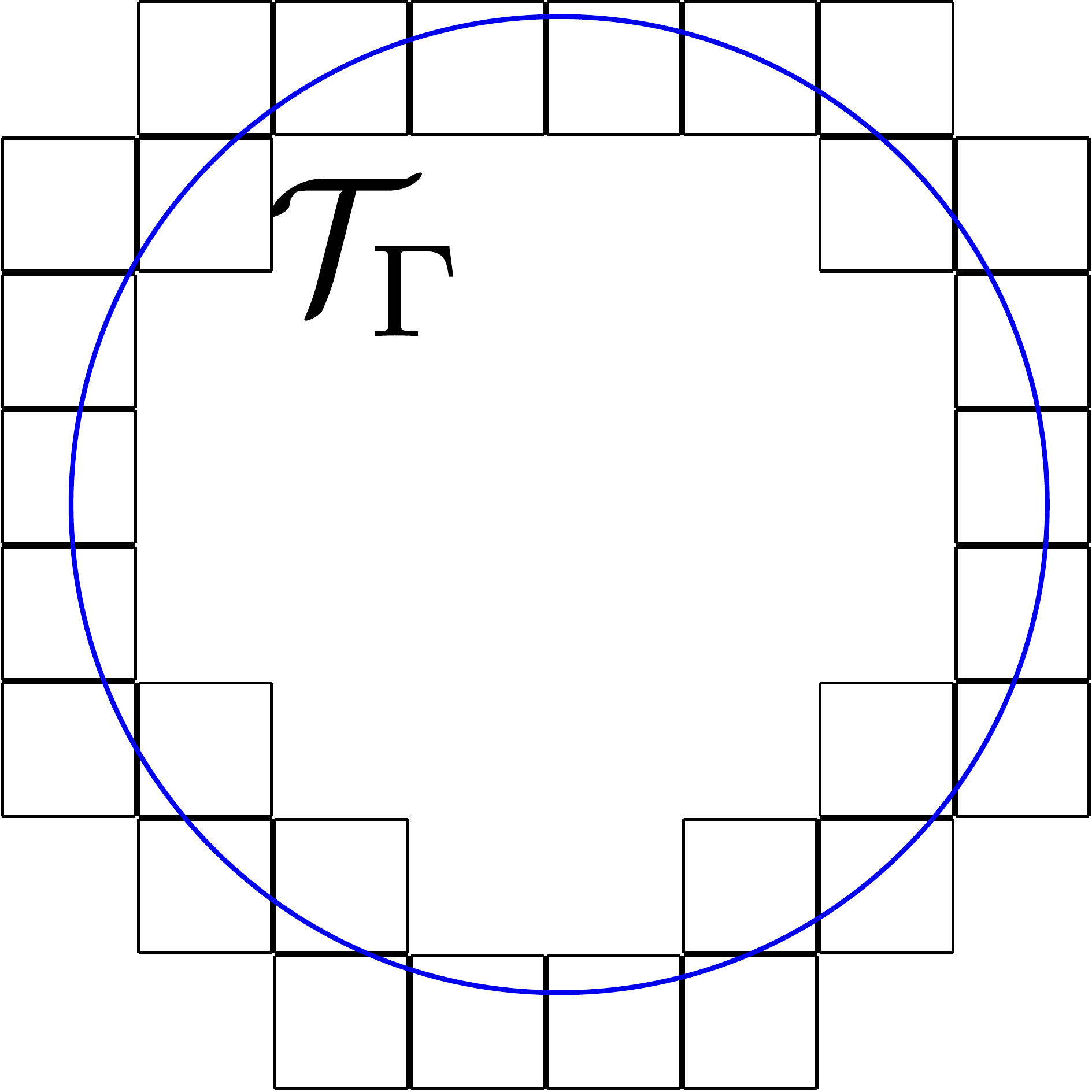}
\caption{\label{fig:cut_mesh}}
\end{subfigure}
\caption{The \textbf{(a)} subdomain $\Omega_1$ immersed in a mesh $\T_{1}$, \textbf{(b)} subdomain $\Omega_2$ immersed in a mesh $\T_{2}$ \siyang{}{and domain $\Omega$ immersed in $\T$}{Comment 4, Reviewer 1, also add $\T$ in figure b)}, and \textbf{(c)} intersected elements $\T_\Gamma$.}
\label{fig:cutFEM_meshes}
\end{figure}

Let $\Omega_s$ be covered by a structured triangulation, $\T_{s}$, so that each element $T\in \T_{s}$ has some part inside of $\Omega_s$; see Figures~\ref{fig:outerTriangulation} and~\ref{fig:innerTriangulation}.
Here, $s = 1,2$ is an index for the composite domain problem (\ref{eqn:PDECD1}--\ref{eqn:PDECD7}), which will be omitted when referring to the single domain problem (\ref{eqn:PDESD1},~\ref{eqn:PDESD2-3}).
(For the latter, note that $\T$ covers $\Omega$.)
\siyang{}{Typically $\T_1$ and $\T_2$ would be created from a larger mesh by removing some of the cells.}{Comment 4, Reviewer 1}
Further, let $\T_\Gamma=\{T\in\T : T \cap \Gamma \neq \emptyset\}$ be the set of intersected elements; see Figure~\ref{fig:cut_mesh}.
In the following, we shall use $\Gamma$ both for the immersed boundary of the single domain problem and for the immersed interface of the composite domain problem, in order to make the connection to the set $\T_\Gamma$ clearer.

To construct the finite element spaces we use Lagrange elements with Gauss--Lobatto nodes of order $p$ ($Q_p$-elements).
Let $V_h^s$ denote a continuous finite element space on $\Omega_s$, consisting of $Q_p$-elements on the mesh $\T_s$:
\begin{equation}
V_h^s=\left \{ v\in C^0(\Omega_s): \evaluated{v}{T}\in Q_p(T), \, T\in \T_s\right \}.
\label{eq:fullFEMspace}
\end{equation}
For the single domain problem (\ref{eqn:PDESD1},~\ref{eqn:PDESD2-3}) we solve for the solution $u \in V_h$; while for the composite domain problem (\ref{eqn:PDECD1}--\ref{eqn:PDECD7}), we solve for the pair $\pair{u}\in V_h^1 \times V_h^2$.
For the latter problem, this means that the degrees of freedom are doubled over elements belonging to $\T_\Gamma$.

%
%
We begin by stating the weak formulation for the single domain problem (\ref{eqn:PDESD1},~\ref{eqn:PDESD2-3}).
Let $(\cdot,\cdot)_X$ and $\bsprod{\cdot}{\cdot}_Y$ be the $L_2$ scalar products taken over the two- and one-dimensional domains $X\subset \R^2$ and $Y\subset \R^{1}$, respectively.
The present method is based on modifying the weak formulation by using Nitsche's method \cite{Nitsche1971} to enforce the boundary condition \eqref{eqn:PDESD2-3}. 
By multiplying \eqref{eqn:PDESD1} with a test function $v\in V_h$, and integrating by parts, we obtain: 
\begin{equation}
(\dot{u},v)_{\Omega}+(\lambda \nabla u, \nabla v)_{\Omega} -\bsprod{\lambda \pdd{u}{n}}{v}_{\Gamma}
=(f,v)_{\Omega}, \quad \forall v\in V_h.
\label{eq:naive_weak_form}
\end{equation}
Note that \eqref{eqn:PDESD2-3} is consistent with the following terms:
\begin{align}
\frac{\gamma_{D}}{\elementSize}\bsprod{\lambda u}{v}_{\Gamma} &= \frac{\gamma_{D}}{\elementSize}\bsprod{\lambda \psi}{v}_{\Gamma},\label{eq:dirichletConsistent1} \\ 
-\bsprod{u}{\lambda \frac{\partial v}{\partial n}}_{\Gamma} &= -\bsprod{\psi}{\lambda \frac{\partial v}{\partial n}}_{\Gamma},\label{eq:dirichletConsistent2}
\end{align}
where $\gamma_D$ is a constant, and $\elementSize$ is the side length of the quadrilaterals in the triangulation. Now, adding (\ref{eq:dirichletConsistent1},~\ref{eq:dirichletConsistent2}) to \eqref{eq:naive_weak_form} gives the following weak form: Find $u \in V_h$ such that 
\begin{equation}
(\dot{u},v)_{\Omega}+a(u,v)=L(v), \quad \forall v\in V_h,
\label{eq:weak_form_short}
\end{equation}
where
\begin{align}
a(u,v) & =(\lambda \nabla u, \nabla v)_{\Omega} -\bsprod{\lambda \pdd{u}{n}}{v}_{\Gamma}
-\bsprod{u}{\lambda \pdd{v}{n}}_{\Gamma}+\frac{\gamma_{D}}{\elementSize} \bsprod{\lambda u}{v}_{\Gamma}, \\
L(v) & =(f,v)_{\Omega}
+\bsprod{\lambda \psi}{\frac{\gamma_{D}}{\elementSize} v- \pdd{v}{n}}_{\Gamma}.
\end{align}

\begin{figure}[tb] 
\centering
\begin{subfigure}[b]{0.3\linewidth}
\includegraphics[width=\columnwidth]{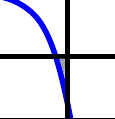}
\caption{\label{fig:small_cut}}
\end{subfigure}
\hfill
\begin{subfigure}[b]{0.3\linewidth}
\includegraphics[width=\columnwidth]{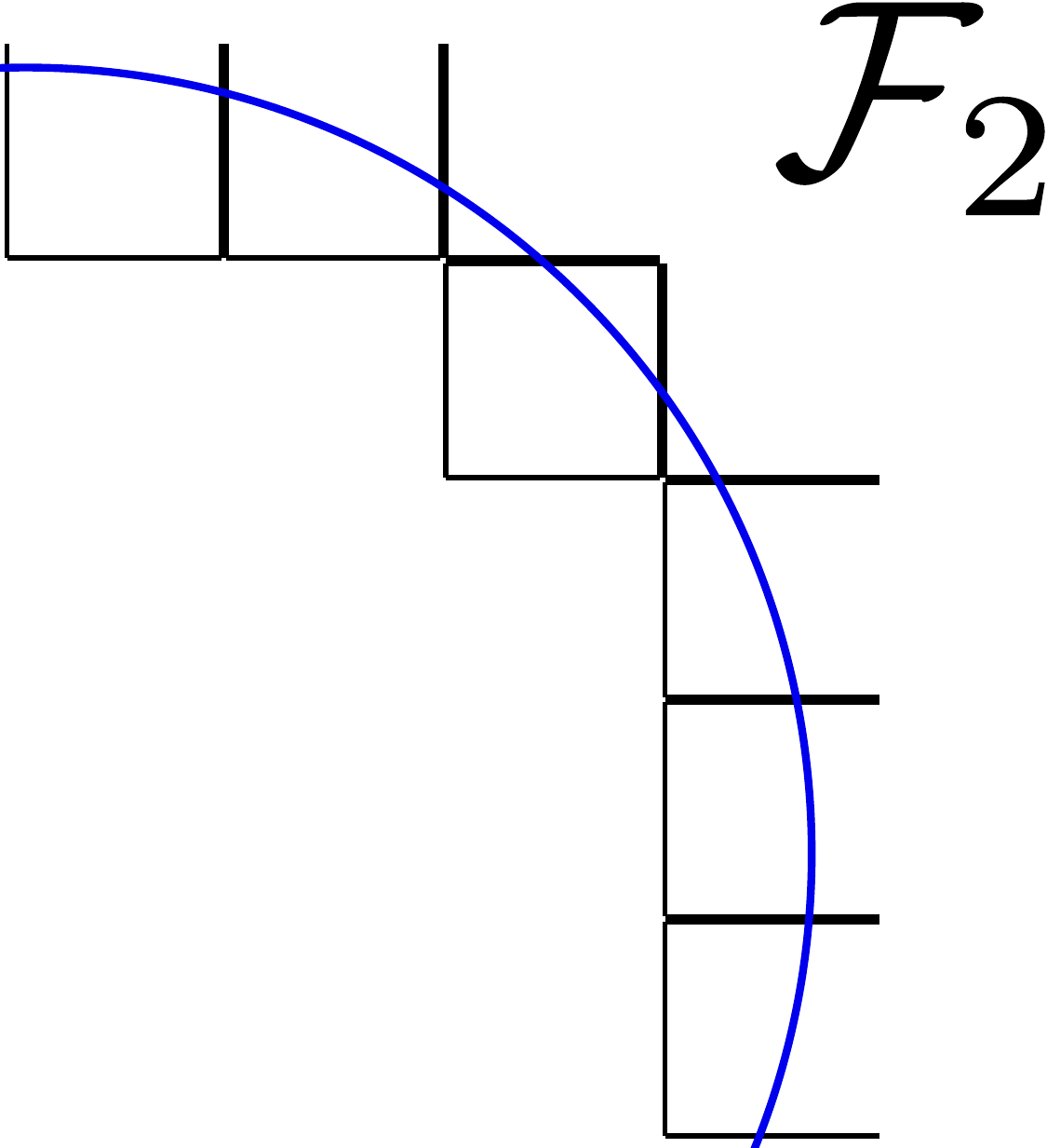}
\caption{\label{fig:inner_stabilized_faces}}
\end{subfigure}
\hfill
\begin{subfigure}[b]{0.3\linewidth}
\includegraphics[width=\columnwidth]{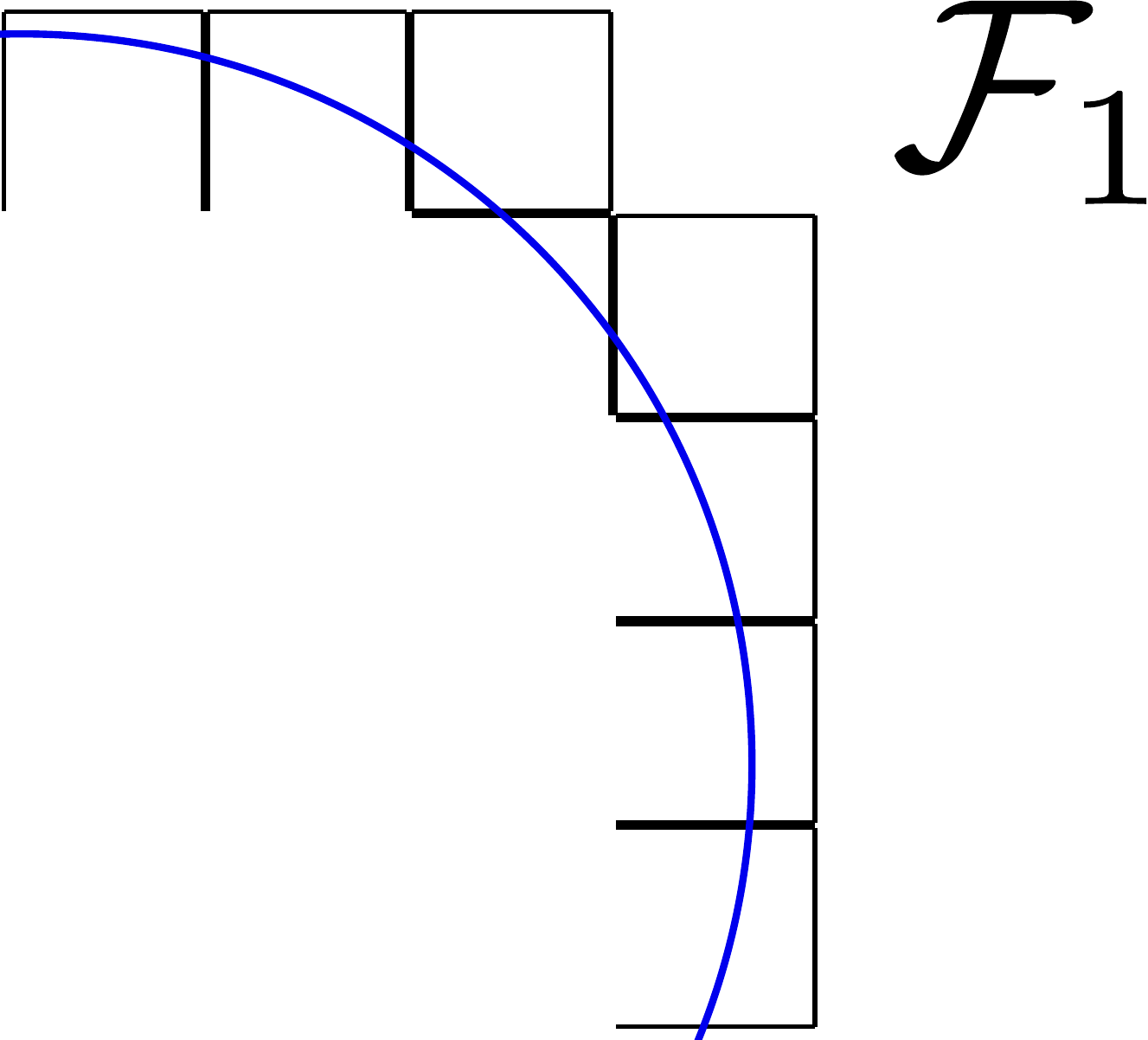}
\caption{\label{fig:outer_stabilized_faces}}
\end{subfigure}
\caption{\textbf{(a)} An element having a small intersection, shown in gray, with the domain; \textbf{(b)} faces belonging to $\Fset_{2}$; and \textbf{(c)} faces belonging to $\Fset_{1}$.}
\label{fig:cut_figs}
\end{figure}

For $\T_\Gamma$ (the elements intersected by $\Gamma$), note that one must integrate only over the part of the element that lies inside $\Omega$.
A problem with this is that one cannot control how the intersections (cuts) between $\Omega$ and $\T$ are made.
Depending on how $\Omega$ is located with respect to the triangulation, some elements can have an arbitrarily small intersection with the domain -- see, for example, Figure~\ref{fig:small_cut}.
If $\Omega$ is moved with respect to $\T$ to make the cut arbitrarily small, then the condition numbers of the mass and stiffness matrices can become arbitrarily large. 

To mitigate this issue, in this work we add a stabilizing term $j$ -- defined shortly in \eqref{eq:jump_stabilization} -- to the mass and stiffness matrices, so that their condition numbers are 
bounded, independently of how the domain $\Omega$ is located with respect to the triangulation $\T$ \cite{burman_hansbo_cut_II,massing_stabilized_2014}. 
Adding stabilization to \eqref{eq:weak_form_short} results in the following weak form: Find $u \in V_h$ such that 
\begin{equation}
(\dot{u},v)_{\Omega}+ \gamma_M j(\dot{u},v)+a(u,v) + \gamma_A \elementSize^{-2}\lambda j(u,v)=L(v), \quad \forall v\in V_h,
\label{eq:weak_form_stabilized}
\end{equation}
where $\gamma_M$  and $\gamma_A$ are scalar constants.

%
%
In order to state the definition of stabilization \eqref{eq:jump_stabilization}, denote by $\Fset_{s}$ the set of faces, as seen in Figures~\ref{fig:inner_stabilized_faces} and~\ref{fig:outer_stabilized_faces}. 
That is, $\Fset_{s}$ is the set of all faces of the elements in $\T_\Gamma$, excluding the boundary faces of $\T_s$: 
\begin{equation}
\Fset_{s}=\{F=T_A \cap T_B \; : \; T_A \in \T_\Gamma \;\text{ or }\; T_B \in \T_\Gamma, \quad T_A,T_B\in \T_{s}  \}. 
\end{equation}
Then, the stabilization term is defined as: 
\begin{equation}
j_s(u,v)=\sum_{F\in \Fset_{s}} 
\sum_{k=1}^p \frac{\elementSize^{2k+1}}{(2k+1)(k!)^2}\bsprod{[\partial_n^k u]}{[\partial_n^k v]}_F,
\label{eq:jump_stabilization}
\end{equation}
where $[u] = u|_{F_+}-u|_{F_-}$  is the jump over a face, $F$; $n$ refers to a normal of $F$;  and $\partial_n^k u$ denotes the $k$-th order normal derivative. 
\siyang{}{
The scaling with respect to $k$ of the terms in \eqref{eq:jump_stabilization} is based on how the stabilization was derived. 
In particular, the $k!$-factors come from the Taylor-expansion and the factor $2k+1$ comes from integrating each term once.}{
Addition to Comment 5, Reviewer 1}

%
%
We now consider the composite domain problem (\ref{eqn:PDECD1}--\ref{eqn:PDECD7}). 
To derive the weak formulation, one follows essentially the same steps as for the single domain problem, namely:
\begin{enumerate}
  \item For both \eqref{eqn:PDECD1} and \eqref{eqn:PDECD2}, multiply the equation for $\idx{u}{s}$ with a test function  $\idx{v}{s} \in V_h^s$, and then integrate by parts; 
  \item Add terms consistent with the interface and boundary conditions; and 
  \item Add stabilization terms $j_1$ and $j_2$ over $\Fset_1$ and $\Fset_2$, respectively.
\end{enumerate}
This results in the following weak formulation for (\ref{eqn:PDECD1}--\ref{eqn:PDECD7}). Find $u=\pair{u} \in V_h^1 \times V_h^2 $ such that:
\begin{align}
M(\dot{u},v)&+A(u,v)+a_\Gamma(u,v)+a_{\partial \Omega}(u,v) \nonumber \\
&=L_\Omega(v)+L_\Gamma(v)+L_{\partial \Omega}, \quad 
\forall v=\pair{v}\in V_h^1 \times V_h^2,
\label{eq:weak_short_form}
\intertext{where the bilinear forms $M$ and $A$ correspond to the stabilized mass and stiffness matrices:}
M(\dot{u},v)&=\sum_{s=1}^2(\idx{\dot{u}}{s},\idx{v}{s})_{\Omega_s} + \gamma_M j_s(\idx{\dot{u}}{s},\idx{v}{s}), \label{eq:weak_sf_2} \\
A(u,v)&=\sum _{s=1}^2(\lambda \nabla \idx{u}{s}, \nabla \idx{v}{s})_{\Omega_s} +\gamma_A \elementSize^{-2} \lambda j_s(\idx{u}{s},\idx{v}{s}); \label{eq:weak_sf_3}
\intertext{$L_\Omega$ corresponds to the forcing function:}
L_\Omega(v) & =\sum _{s=1}^2(f_s,\idx{v}{s})_{\Omega_s}; \label{eq:weak_sf_4} \\
\intertext{$a_\Gamma$ and $L_\Gamma$ consistently enforce the interface conditions (\ref{eqn:PDECD6},~\ref{eqn:PDECD7}):}
a_\Gamma (u,v) & =-\bsprod{[u]}{ \{\lambda \pdd{v}{n}\}}_{\Gamma}-\bsprod{ \{\lambda \pdd{u}{n}\}}{[v]}_{\Gamma} + \bsprod{\frac{\gamma_\Gamma}{\elementSize} [u]}{[v]}_{\Gamma}, \label{eq:weak_sf_5} \\
L_\Gamma(v) & =\bsprod{\frac{\gamma_\Gamma}{\elementSize} \mu_1}{[v]}_\Gamma + \bsprod{\kappa_1\mu_2}{\idx{v}{2}}_\Gamma + \bsprod{\kappa_2\mu_2}{\idx{v}{1}}_\Gamma -\bsprod{\mu_1}{\{\lambda \pdd{v}{n}\}}_\Gamma; \label{eq:weak_sf_6}
\intertext{and the terms $a_{\partial \Omega}$ and $L_{\partial \Omega}$ enforce the boundary condition \eqref{eqn:PDECD5} along the outer boundary, $\partial \Omega$:}
a_{\partial \Omega}(u,v) & =-\bsprod{\lambda \pdd{\idx{u}{1}}{n}}{\idx{v}{1}}_{\partial \Omega} -\bsprod{\idx{u}{1}}{\lambda \pdd{\idx{v}{1}}{n}}_{\partial \Omega}+\frac{\gamma_{D}}{\elementSize} \bsprod{\lambda \idx{u}{1}}{\idx{v}{1}}_{\partial \Omega}, \label{eq:weak_sf_7} \\
L_{\partial \Omega}(v) & =\bsprod{\lambda\psi}{\frac{\gamma_D}{\elementSize} \idx{v}{1} - \pdd{\idx{v}{1}}{n}}_{\partial \Omega}. \label{eq:weak_sf_8}
\end{align}
In (\ref{eq:weak_sf_5}--\ref{eq:weak_sf_8}), $n$ denotes the outward pointing normal at either $\Gamma$ or $\partial\Omega$ (depending on the domain of integration); $\kappa_1 + \kappa_2 = 1$, so that $\{ v \} = \kappa_1 v_1+ \kappa_2 v_2$ is a convex combination; and
$\gamma_\Gamma$, $\kappa_1$, $\kappa_2$ are chosen as in \cite{cutfem_2015}:
\begin{equation}
\kappa_1 = \frac{\idx{\lambda}{2}}{\idx{\lambda}{1} + \idx{\lambda}{2}}, \quad 
\kappa_2 = \frac{\idx{\lambda}{1}}{\idx{\lambda}{1} + \idx{\lambda}{2}}, \quad
\gamma_{\Gamma} = \gamma_D \frac{\idx{\lambda}{1}\idx{\lambda}{2}}{\idx{\lambda}{1} + \idx{\lambda}{2}}. \label{eq:weak_sf_params1}
\end{equation}
The remaining parameters (appearing in Equations~\ref{eq:weak_sf_2},~\ref{eq:weak_sf_3},~\ref{eq:weak_sf_7}--\ref{eq:weak_sf_params1}) are given by:
\begin{equation}
\gamma_M=0.75,\quad \gamma_A=1.5, \quad \gamma_D=5p^2.
\end{equation}
The scaling of $\gamma_D$ with respect to $p$ follows from an inverse inequality. 
\siyang{and the constants are based on numerical experiments,}{
When $p=1$ these reduce to the same parameters as the ones used in \cite{stickoLowerOrder2016}, where $\gamma_M$ was chosen based on numerical experiments on the condition number of the mass matrix.
This also agrees with the choice of $\gamma_A$ and $\gamma_D$ in \cite{burman_hansbo_cut_II}, where $\gamma_A$ was investigated numerically.}{
Comment 5, Reviewer 1}

In order to use cut--FEM, one needs a way to perform integration over the intersected elements $\T_\Gamma$.
For example, with the interface problem, on each element $K \in \T_\Gamma$, we need a quadrature rule for the $K\cap \Omega_1$, $K\cap \Omega_2$ and $K\cap \Gamma$.
For the numerical tests in this work (Section~\ref{sect:TestProblems}), we represent the geometry by a level set function, and compute high-order accurate quadrature rules with the algorithm from \cite{saye_highOrder_2015}.

\begin{remark}{}
Optimal (second-order) convergence was rigorously proven for cut--FEM applied to the Poisson problem in \cite{burman_hansbo_cut_II}.
As far as we know, there is no rigorous proof of higher-order convergence for cut--FEM, though such a proof would likely be similar to the second-order case.
\end{remark}

  \subsection{DPM} \label{sec:Overview:DPM}



\makeatletter
\renewcommand*\env@matrix[1][\arraystretch]{%
  \edef\arraystretch{#1}%
  \hskip -\arraycolsep
  \let\@ifnextchar\new@ifnextchar
  \array{*\c@MaxMatrixCols c}}
\makeatother


We continue in this section with a brief introduction to the Difference Potentials Method (DPM),  
{which was originally proposed by V. S. Ryaben'kii (see \cite{Ryabenkii1988,Ryabenkii2002,Ryabenkii2012}, and see \cite{EpshteynSofronovTsynkov2015,Godunov2015} for papers in his honor).} 
Our aim is to consider the numerical approximation of PDEs on arbitrary, smooth geometries (defined either explicitly or implicitly) using {the DPM} together with standard, finite-difference discretizations of (\ref{eqn:PDESD1}) or (\ref{eqn:PDECD1},~\ref{eqn:PDECD2}) on uniform, Cartesian grids, which need not conform with boundaries or interfaces. 
To this end, we work with high-order methods for interface problems based on Difference Potentials, which were originally developed in \cite{ryaben2006algorithm} and \cite{phdthesisJA,albright2017high,AlbrightEpshteynSteffen2015,albright2015high,Epshteyn2014,epshteyn2015high}.  We also introduce new developments here for handling implicitly-defined geometries.  
(The reader can consult \cite{Ryabenkii2002} for the general theory of the Difference Potentials Method.)

Broadly, the main idea of the DPM is to reduce uniquely solvable and well-posed boundary value problems {in a domain $\Omega$} to pseudo-differential \textit{Boundary Equations with Projections} (BEP) {on the boundary of $\Omega$}.  
{First, we introduce a computationally simple auxiliary domain as part of the method.  The original domain is embedded into the auxiliary domain, which is then discretized using a uniform Cartesian grid.} 
{Next, we define a Difference Potentials operator \textit{via} the solution of a simple Auxiliary Problem (defined on the auxilairy domain), and construct the discrete, pseudo-differential \textit{Boundary Equations with Projections} (BEP) at grid points near the continuous boundary or interface $\Gamma$. 
(This set of grid points is called the \textit{discrete grid boundary}.)} 
Once constructed, the BEP are then solved 
{together with the boundary/interface conditions} 
to obtain the value of the solution at 
{the discrete grid boundary.} 
Lastly, using these 
{reconstructed values of the solution at the discrete grid boundary}, 
the approximation to the solution in the domain $\Omega$ is 
{obtained through} 
the \textit{discrete, generalized Green's formula}. 

Mathematically, the DPM is a discrete analog of the method of Calder\'{o}n's potentials in the theory of partial differential equations. 
The DPM, however, does not require 
explicit knowledge of Green's functions. 
Although we use an Auxiliary Problem (AP) discretized by finite differences, the DPM is not limited to this choice of spatial discretization. 
Indeed, numerical methods based on the idea of Difference Potentials can be designed with whichever choice of spatial discretization is most natural for the problem at hand (\textit{e.g.}, see \cite{Epshteyn2012}). 

Practically, the main computational complexity of the DPM reduces to the required solutions of the AP, which can be done very efficiently using fast, standard $\mathcal{O}(N \log N)$ solvers. 
Moreover, in general the DPM can be applied to problems with general boundary or interfaces conditions, with no change to the discretization of the PDE. 

Let us now briefly introduce the DPM for the numerical approximation of parabolic interface models (\ref{eqn:PDECD1}--\ref{eqn:PDECD7}).  First, we must introduce the point-sets that will be used throughout the DPM.  
(Note that the main construction of the method below applies to the single domain problem (\ref{eqn:PDESD1},~\ref{eqn:PDESD2-3}), after omitting the index $s$ and replacing interface conditions with boundary conditions; see \cite{albright2015high}.) 

Let $\Omega_s$ ($s=1,2$) be embedded in a rectangular auxiliary domain $\Omega_{s}^0$. 
{Introduce a uniform, Cartesian grid denoted $M_s^0$ on $\Omega_s^0$, with grid-spacing $h_s$.} 
{Let $M_s^+ = M_s^0 \cap \Omega_s$ denote the grid points inside each subdomain $\Omega_s$, and} 
{$M_s^- = M_s^0 \setminus M_s^+$ the grid points outside each subdomain $\Omega_s$.} 
Note that the auxiliary domains $\Omega_1^0$, $\Omega_2^0$ 
{and auxiliary grids $M_1^0$, $M_2^0$} 
need not agree, and indeed may be selected completely independently, given considerations regarding accuracy, adaptivity, or efficiency. 

{Define a finite-difference stencil $N_{j, k}^{s, \alpha}$, with $\alpha = 5, 9$, to be the stencil of the standard five-point or a wide nine-point Laplacian, \textit{i.e.},}  
%
\begin{equation} \label{eqn:DPM:FiniteDifferenceStencil}
\begin{aligned}
    N_{j, k}^{s, 5} &= \{ (x_j,y_k), (x_{j\pm1},y_k), (x_j, y_{k\pm1}) \} \quad \text{or} \\
    N_{j, k}^{s, 9} &= \{ (x_j,y_k), (x_{j\pm1},y_k), (x_j, y_{k\pm1}), (x_{j\pm2},y_k), (x_j, y_{k\pm2}) \}. 
\end{aligned}
\end{equation}%
Next, with $\alpha$ fixed, 
define the point-sets %
\begin{equation}
    N_s^0 = \smashoperator{\bigcup_{(x_j, y_k) \in M_s^0}} N_{j, k}^{s, \alpha}, \quad N_s^+ = \smashoperator{\bigcup_{(x_j, y_k) \in M_s^+}} N_{j, k}^{s, \alpha}, \quad \text{and} \quad N_s^- = \smashoperator{\bigcup_{(x_j, y_k) \in M_s^+}} N_{j, k}^{s, \alpha}. 
\end{equation}%
{The point-set $N_s^+$ ($N_s^-$, $N_s^0$) enlarges the point-set $M_s^+$ ($M_s^-$, $M_s^0$), by taking the union of finite-difference stencils at every point in $M_s^+$ ($M_s^-$, $M_s^0$).  Therefore, $N_s^+$ contains a ``thin row'' of points belonging to the complement of $\Omega_s$, so that $N_s^+ \not\subset \Omega_s$, even though $M_s^+ \subset \Omega_s$.  (Likewise, $N_s^- \not\subset \Omega \setminus \Omega_s$, even though $M_s^- \subset \Omega \setminus \Omega_s$). }

\begin{figure}
    \centering
    \begin{subfigure}[b]{0.45\linewidth}
        \centering
        \includegraphics[width=\textwidth,page=1]{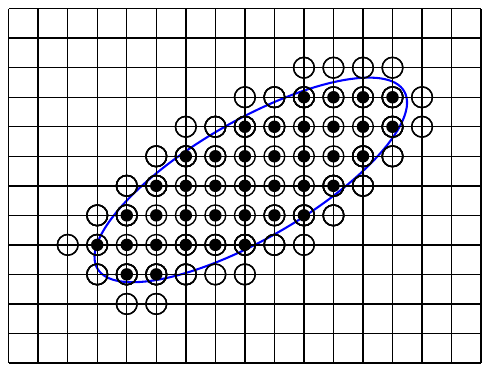}
        \caption{}
        \label{fig:IntroDPM:PointSets:MN}
    \end{subfigure}\quad%
    \begin{subfigure}[b]{0.45\linewidth}
        \centering
        \includegraphics[width=\textwidth,page=1]{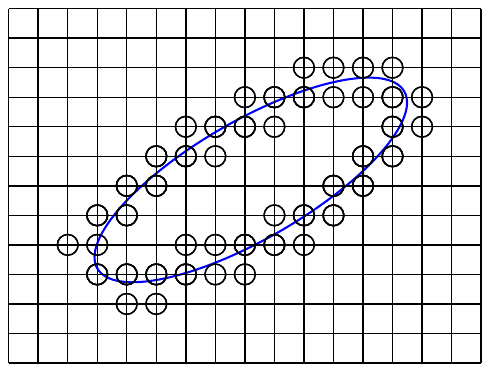}
        \caption{}
        \label{fig:IntroDPM:PointSets:gamma}
    \end{subfigure}

    \caption{An example of the point-sets for the second-order Difference Potentials Method, applied to a single domain ($s$ omitted for brevity), with a rotated ellipse $\Omega$ and a rectangular auxiliary domain $\Omega^0$, showing \textbf{(a)} The physical grids $M^+$ (solid dots) and $N^+$ (open circles), {with $M^+\subset N^+$}, and \textbf{(b)} The discrete grid boundary $\gamma$ (open circles).}
    \label{fig:IntroDPM:PointSets}
\end{figure}

{Lastly, we now define the important point-set %
\begin{equation}
    \gamma_s = N_s^+ \bigcap N_s^-, 
\end{equation}%
which we call the \textit{discrete grid boundary}.  In words, $\gamma_s$ is the set of grid points that straddle the continuous interface $\Gamma$.  (See Figure~\ref{fig:IntroDPM:PointSets} for an example of these points-sets, given a single elliptical domain $\Omega$.)  Note that the point-sets $M_s^+$, $N_s^+$, and $\gamma_s$ will be used throughout the Difference Potentials Method.} 

Here, we define the fully-discrete finite-difference discretization of (\ref{eqn:PDECD1}, \ref{eqn:PDECD2}), and then define the Auxiliary Problem. 
Indeed, the discretization we consider is %
\begin{equation} \label{eqn:IntroDPM:FullyDiscreteCD}
    L_{\Delta t, h}^s u_{s}^{i+1} = F_s^{i+1}, \quad (x_j, y_k) \in M_s^+, 
\end{equation}%
where (i) $L_{\Delta t, h}^s u_{s}^{i+1} := \lambda_s(t^{i+1}) \Delta_h u_{s}^{i+1} - \sigma u_{s}^{i+1}$, (ii) $\Delta_h$ is either a five- or nine-point Laplacian 
{in each subdomain}, 
and (iii) $\sigma$ and $F_s^{i+1}$ follow from the choice of time- and spatial-discretizations.  (Here, we have simplified notation slightly by assuming that $h := h_1 = h_2$, which need not be the same in general.)  {For full details of the discretization, 
{including the choice of BDF2 or BDF4 in the time discretization,} 
we refer the reader to Appendix~\ref{appendix:DPM:discretePDE}.} 

The choices of discretization \eqref{eqn:IntroDPM:FullyDiscreteCD} in each subdomain need not be the same.  As in \cite{phdthesisJA,albright2015high}, one could choose a second- and fourth-order discretization on $M_1^+$ and $M_2^+$, respectively, given considerations about accuracy, adaptivity, expected regularity of the analytical solution in each domain, etc.  

{Next, we define the discrete Auxiliary Problem, which plays a central role in the construction of the 
{Difference Potentials operator, the resulting Boundary Equations with Projection at the discrete grid boundary,}
and in the numerical approximation of the solution \textit{via} the discrete, generalized Green's formula.} 

\begin{definition}[Discrete Auxiliary Problem (AP)] \label{def:IntroDPM:AP}
    At time $t^{i+1}$, given the right-hand side grid function $q_s^{i+1}$: $M_s^0 \to \mathbb{R}$, the following difference equations (\ref{eqn:IntroDPM:AP1},~\ref{eqn:IntroDPM:AP2}) are defined as the discrete AP. 
    \begin{alignat}{2}
        L_{\Delta t, h}^s u_s^{i+1} &= q_s^{i+1}, && \quad (x_j,y_k)\in M^0_s \label{eqn:IntroDPM:AP1} \\
        u_s^{i+1} &=0, && \quad (x_j,y_k)\in N^0_s\backslash M^0_s \label{eqn:IntroDPM:AP2} 
    \end{alignat}
\end{definition}

\begin{remark} \label{rmk:IntroDPM:AP}
    {For a given right-hand side $q_s^{i + 1}$, the solution of the discrete AP (\ref{eqn:IntroDPM:AP1},~\ref{eqn:IntroDPM:AP2}) defines a discrete Green's operator $G_{\Delta t, h}^s q_s^{i + 1}$. 
    The choice of boundary conditions \eqref{eqn:IntroDPM:AP2} will affect the resulting grid function $G_{\Delta t, h}^s q_s^{i + 1}$, and thus the Boundary Equations with Projection defined below. 
    However, the choice of boundary conditions 
    {\eqref{eqn:IntroDPM:AP2} in the AP} 
    will not affect the numerical approximation of (\ref{eqn:PDECD1}--\ref{eqn:PDECD7}), so long as the discrete AP is uniquely solvable and well-posed.}
\end{remark}

Let us denote by $G_{\Delta t, h}^s F_s^{i + 1}$ the particular solution on $N^+_s$ of the fully-discrete problem \eqref{eqn:IntroDPM:FullyDiscreteCD}, defined by solving the AP (\ref{eqn:IntroDPM:AP1},~\ref{eqn:IntroDPM:AP2}) with %
\begin{equation} \label{eqn:DPM:rhs_PS}
    q_s^{i + 1} = \begin{cases}
        F_s^{i + 1}, & (x_j, y_k) \in M_s^+, \\
        0, & (x_j, y_k) \in M_s^-, 
    \end{cases}
\end{equation}%
and restricting the solution from $N_s^0$ to $N_s^+$.

Let us also introduce a linear space $\mathbf{V}_{\gamma_s}$ of all grid functions denoted $v_{\gamma_s}^{i + 1}$, which are defined on $\gamma_s$ and extended by zero to the other points of $N_s^0$.  These grid functions are referred to as \textit{discrete densities on $\gamma_s$}.

{\begin{definition}[The Difference Potential of a density.]
    The Difference Potential of a given density $v_{\gamma_s}^{i + 1}$ is the grid function $P_{N_s^+}^{i + 1} v_{\gamma_s}^{i + 1}$ on $N_s^+$, defined by solving the AP (\ref{eqn:IntroDPM:AP1},~\ref{eqn:IntroDPM:AP2}) with %
    \begin{equation} \label{eqn:DPM:rhs_DP}
        q_s^{i + 1} = \begin{cases}
            0, & (x_j, y_k) \in M_s^-, \\
            L_{\Delta t, h}^s [v_{\gamma_s}^{i + 1}], & (x_j, y_k) \in M_s^+, 
        \end{cases}
    \end{equation}%
    and restricting the solution from $N_s^0$ to the point-set $N_s^+$. 
\end{definition}}

{Note that $P_{N_s^+}^{i + 1} : \mathbf{V}_{\gamma_s} \to N_s^+$ is a linear operator on the space $\mathbf{V}_{\gamma_s}$ of densities. 
Moreover, the coefficients of $P_{N_s^+}^{i + 1}$ can be computed by solving the AP (Definition~\ref{def:IntroDPM:AP}) with the appropriate density $v_{\gamma_s}^{i+1}$ defined at the points $(x_{j}, y_{k}) \in \gamma_s$.} 

{\begin{definition}[The Trace operator.] \label{def:DPM:TraceOp}
    Given a grid function $v_s^{i + 1}$, we denote by $\Tr_{\gamma_s} [v_s^{i + 1}]$ the Trace (or Restriction) from $N_s^+$ to $\gamma_s$. 

    Moreover, for a given density $v_{\gamma_s}^{i + 1}$, denote the trace of the Difference Potential of $v_{\gamma_s}^{i + 1}$ by $P_{\gamma_s} v_{\gamma_s}^{i + 1}$.  In other words, $P_{\gamma_s} v_{\gamma_s}^{i + 1} = \Tr_{\gamma_s} [P_{N_s^+}^{i + 1} v_{\gamma_s}^{i + 1}]$. 
\end{definition}}

{Now we can state the central theorem of the Difference Potentials Method that will allow us to reformulate the finite-difference equations \eqref{eqn:IntroDPM:FullyDiscreteCD} on $M_s^+$ (without imposing any boundary or interface conditions yet) into the equivalent Boundary Equations with Projections on $\gamma_s$.} 

\begin{theorem}[Boundary Equations with Projection (BEP)] \label{thm:IntroDPM:BEPEquivalence}
At time-level $t^{i+1}$, the discrete density $u^{i+1}_{\gamma_s}$ ($s=1,2$) is the trace of some solution $u^{i+1}_s$ on domain $\Omega_s$ to the Difference Equations~\eqref{eqn:IntroDPM:FullyDiscreteCD}, \textnormal{i.e.}, $u^{i+1}_{\gamma_s}:={\Tr_{\gamma_s}}[u^{i+1}_s]$, if and only if the following BEP holds
\begin{align}\label{eqn:IntroDPM:BEP}
u^{i+1}_{\gamma_s}-P_{\gamma_s}^{i+1}u^{i+1}_{\gamma_s}={\Tr_{\gamma_s}}[G^{i+1}_{\Delta t, h}F^{i+1}_s], \quad (x_j,y_k)\in\gamma_s,
\end{align}
with $\Tr_{\gamma_s}[\cdot]$ and $P_{\gamma_s}$ defined in Definition~\ref{def:DPM:TraceOp}. 
\end{theorem}

\begin{proof}
    See \cite{Ryabenkii2002} for the general theory of DPM (including the proof for general elliptic PDE), or one of \cite{phdthesisJA,AlbrightEpshteynSteffen2015,albright2015high} for the proof in the case of parabolic interface problems. 
\end{proof}

{\begin{remark}
    A given density $v_{\gamma_s}^{i+1}$ is the trace of some solution of the fully-discrete finite-difference equations \eqref{eqn:IntroDPM:FullyDiscreteCD} if and only if it is a solution of the BEP. 

    However, since boundary or interface conditions have not yet been imposed, the BEP will have infinitely many solutions $u_{\gamma_s}^{i+1}$.  As originally disucssed in \cite{phdthesisJA,albright2017high,AlbrightEpshteynSteffen2015,albright2015high,Epshteyn2014,epshteyn2015high}, in this work we consider the following approach in order to find a unique solution of the BEP. 

    At each time level $t^{i + 1}$, one can approximate the solution of (\ref{eqn:PDECD1}--\ref{eqn:PDECD7}) at the discrete grid boundary $\gamma_s$, using the Cauchy data of (\ref{eqn:PDECD1}--\ref{eqn:PDECD7}) on the continuous interface $\Gamma$, up to the desired second- or fourth-order accuracy.  (By Cauchy data, 
    we mean the trace of the solution of (\ref{eqn:PDECD1}--\ref{eqn:PDECD7}), together with the trace of its normal derivative, on $\Gamma$.) Below, we will define an Extension Operator 
    which will extend the Cauchy data of (\ref{eqn:PDECD1}--\ref{eqn:PDECD7}) from $\Gamma$ to $\gamma_s$. 

    As we will see, the Extension Operator in this work depends only on the given parabolic interface model.  Moreover, we will use a finite-dimensional, spectral representation for the Cauchy data of (\ref{eqn:PDECD1}--\ref{eqn:PDECD7}) on $\Gamma$.  
    {Then, we will use the Extension Operator, together with the BEP \eqref{eqn:IntroDPM:BEP} and the interface conditions (\ref{eqn:PDECD6},~\ref{eqn:PDECD7}), to obtain a linear system of equations for the coefficients of the finite-dimensional, spectral representation.  Hence, the derived BEP will be solved for the unknown coefficients of the Cauchy data.  Using this obtained Cauchy data, we will construct the approximation of (\ref{eqn:PDECD1}--\ref{eqn:PDECD7}) using the Extension Operator, together with the discrete, generalized Green's formula.} 
\end{remark}}



%
Let us now briefly discuss the Extension Operator for the second-order numerical method, and refer the reader to Appendix~\ref{appendix:DPM:Extop} for {details (including details for the fourth-order numerical method)}.  For points in the vicinity of $\Gamma$, we define a coordinate system $(d, \vartheta)$, where $\vartheta$ is arclength from some reference point, and {$d$ is the signed distance in the normal direction from the point to $\Gamma$.} 
Now, as a first step towards defining the Extension Operator, we define a new function %
\begin{equation} \label{eqn:ExtOperatorStep1}
    v_s^{i+1}(d, \vartheta) = v_s^{i+1}(0, \vartheta) + \sum_{l = 1}^p \frac{1}{l!} \frac{\partial^l v_s^{i+1}(0, \vartheta)}{\partial n^l} d^l, 
\end{equation}%
{where $n$ is the unit outward normal vector at $\Gamma$.} 
{We choose $p = 2$ for the second-order method (which we will discuss now) and $p = 4$ for the fourth-order method (see Appendix~\ref{appendix:DPM:Extop}).} 

As a next step for the second-order method (BDF2--DPM2), we define %
\begin{equation} \label{eqn:ExtOperatorStep2}
    v_s^{i+1}(0, \vartheta) = u_s^{i + 1}\rvert_\Gamma, \quad \frac{\partial v_s^{i+1}(0, \vartheta)}{\partial n} = \left. \frac{\partial u_s^{i + 1}}{\partial n} \right|_\Gamma, \quad 
    \text{and} \quad \frac{\partial^2 v_s^{i+1}(0, \vartheta)}{\partial n^2} = \left. \frac{\partial^2 u_s^{i + 1}}{\partial n^2} \right|_\Gamma, 
\end{equation}%
where $u_s^{i + 1} := u_s(x, y, t^{i + 1})$, $\frac{\partial u_s^{i + 1}}{\partial n} := \frac{\partial u_s(x, y, t^{i + 1})}{\partial n}$, etc. 
As a last step, a straightforward sequence of calculations (see Appendix~\ref{appendix:DPM:Extop}) shows that %
\begin{equation} \label{eqn:ExtOperatorStep3}
        \frac{\partial^2 u_s^{i + 1}}{\partial n^2} \approx \frac{1}{\lambda_s} \left( \frac{3 u_s^{i + 1} - 4 u_s^i + u_s^{i - 1}}{2 \Delta t} - f^{i + 1} \right) - \frac{\partial^2 u_s^{i + 1}}{\partial \vartheta^2} + \kappa \frac{\partial u_s^{i+1}}{\partial n},
\end{equation}%
{where $\kappa$ denotes the curvature of $\Gamma$.} 
Therefore, with $v_s^{i+1}(d, \vartheta)$ defined by (\ref{eqn:ExtOperatorStep1}--\ref{eqn:ExtOperatorStep3}), the only unknown data at each time step $t^{i + 1}$ are the unknown Dirichlet data $u_s^{i + 1}$ and the unknown Neumann data $\frac{\partial u_s^{i + 1}}{\partial n}$. 
The Extension Operator will incorporate the interface conditions (\ref{eqn:PDECD6},~\ref{eqn:PDECD7}) when it is combined with the BEP \eqref{eqn:IntroDPM:BEP}, so that the only independent unknowns at each time step $t^{i + 1}$ will be $\left. \left(u_1^{i + 1}, \frac{\partial u_1^{i + 1}}{\partial n} \right) \right|_\Gamma$ or $\left. \left(u_2^{i + 1}, \frac{\partial u_2^{i + 1}}{\partial n} \right) \right|_\Gamma$. 
(This is also true for the fourth-order numerical method -- see Appendices~\ref{appendix:DPM:Extop} and~\ref{appendix:DPM:CompleteSystemOfEquations}.) 

Now we are ready to define the Extension Operator that extends the Cauchy data of (\ref{eqn:PDECD1}--\ref{eqn:PDECD7}) from $\Gamma$ to $\gamma_s$. 

%
\begin{definition}[The Extension Operator.] \label{def:DPM:ExtOp}
    Let $v_s^{i + 1}(d, \vartheta)$ be defined by (\ref{eqn:ExtOperatorStep1}--\ref{eqn:ExtOperatorStep3}).  Let $\mathfrak{u}_{s, \Gamma}^{i + 1}$ denote the Cauchy data of (\ref{eqn:PDECD1}--\ref{eqn:PDECD7}) at $t^{i + 1}$ on $\Gamma$, \textit{i.e.}, $\mathfrak{u}_{s, \Gamma}^{i + 1} = \left. \left(u_s^{i + 1}, \frac{\partial u_s^{i + 1}}{\partial n} \right) \right|_\Gamma$.  The extension operator $\operatorname{Ex}_s$ that extends $\mathfrak{u}_{s, \Gamma}^{i + 1}$ from $\Gamma$ to $\gamma_s$ is %
    \begin{equation} \label{eqn:DPM:ExtOp:DiscreteVersion}
        \operatorname{Ex}_s \mathfrak{u}_{s, \Gamma}^{i + 1} := v_s^{i + 1}(d, \vartheta)\rvert_{\gamma_s}. 
    \end{equation}%
    {For a given point $(x_j, y_k) \in \gamma_s$, note that $d$ is the signed distance between $(x_j, y_k)$ and its orthogonal projection on $\Gamma$, while $\vartheta$ is the arclength along $\Gamma$ between a reference point and the orthogonal projection of $(x_j, y_k)$.} 
\end{definition}

%
Next, we briefly discuss the finite-dimensional, spectral representation of Cauchy data $\mathfrak{u}_{s, \Gamma}^{i + 1}$.  Indeed, we wish to choose a basis $\phi_{\nu}(\vartheta)$ on $\Gamma$ ($\nu = 1, 2, 3, \ldots$) in order to accurately approximate the two components of the Cauchy data $\mathfrak{u}_{s, \Gamma}^{i + 1}$.  To be specific, whichever basis we choose, we require that %
\begin{align} 
    \varepsilon_{\mathcal{N}^0, \mathcal{N}^1}(\mathfrak{u}_{s, \Gamma}^{i + 1}) &=\operatorname*{min}_{c_{1, \nu}^{s, i+1}, c_{2, \nu}^{s, i+1}} 
    \int_{\Gamma} \Bigg(\Big| u_s^{i + 1} - \sum_{\nu=1}^{\mathcal{N}^0} c_{1, \nu}^{s, i+1} \phi_{\nu}(\vartheta)\Big|^2 \notag \\
    & \qquad\qquad\qquad\qquad \cdots + \Big| \frac{\partial u_s^{i + 1}}{\partial n} -\sum_{\nu=1}^{\mathcal{N}^1} c_{2, \nu}^{s, i+1} \phi_{\nu}(\vartheta)\Big|^2\Bigg) {\rm d\vartheta} \label{eqn:DPM:ResidualSpectralRepr}
\end{align}%
tends to zero as $\mathcal{N}^0, \mathcal{N}^1 \to \infty$, for some sequence of real numbers $(c_{1, \nu}^{s, i+1})_{\nu=1}^{\mathcal{N}^0}$ and $(c_{2, \nu}^{s, i+1})_{\nu=1}^{\mathcal{N}^1}$.  In other words, we require %
\begin{equation}
    \lim_{\mathcal{N}^0, \mathcal{N}^1 \to \infty} \varepsilon_{\mathcal{N}^0, \mathcal{N}^1}(\mathfrak{u}_{s, \Gamma}^{i + 1}) = 0. 
\end{equation}%

Now let us discuss a choice of basis.  In this work, recall that we consider interfaces $\Gamma$ that 
{are at least $C^2(\Gamma)$} 
(due to the choice of smooth, circular geometries).  
Also, as we will see in Section~\ref{sec:TestProblems:List}, each function $u$ considered in the test problems on a composite domain (\ref{eqn:TP2A},~\ref{eqn:TP2B},~\ref{eqn:TP2C}) is locally smooth, in the sense that $u\rvert_{\Omega_1} = u_1$ and $u\rvert_{\Omega_2} = u_2$ are smooth in $\Omega_1$ and $\Omega_2$, respectively.  Moreover, each component of the Cauchy data $\mathfrak{u}_{1, \Gamma}^{i + 1}$ and $\mathfrak{u}_{2, \Gamma}^{i + 1}$ are smooth, periodic functions of arclength $\vartheta$.  (Note that $\mathfrak{u}_{1, \Gamma}^{i + 1}$ and $\mathfrak{u}_{2, \Gamma}^{i + 1}$ need not agree, and indeed do not -- neither $\mu_1(x, y, t)$ nor $\mu_2(x, y, t)$ in (\ref{eqn:PDECD6},~\ref{eqn:PDECD7}) are identically equal to zero, for any of our test problems on a composite domain.)  
%
%

Therefore, in this work, we choose a standard trigonometric basis $\phi_\nu(\vartheta)$, with %
\begin{equation} \label{eqn:DPM:InterfaceFourierSeries}
    \phi_1(\vartheta) = 1, \quad \phi_{2k}(\vartheta) = \cos\left(\frac{2 \pi k \vartheta}{|\Gamma|}\right), \quad \phi_{2k + 1}(\vartheta) = \sin\left(\frac{2 \pi k \vartheta}{|\Gamma|}\right), 
\end{equation}%
and $k > 1$. 
Moreover, at every time step $t^{i + 1}$, we will discretize the Cauchy data $\mathfrak{u}_{s, \Gamma}^{i + 1} = \left. \left(u_s^{i + 1}, \frac{\partial u_s^{i + 1}}{\partial n} \right) \right|_\Gamma$ using this basis.  Therefore, we let %
\begin{equation}\label{eqn:DPM:DiscretizationOfCauchyData}
    \tilde{\mathfrak{u}}_{s, \Gamma}^{i+1} = \sum_{\nu=1}^{\mathcal{N}^0} c^{0, i+1}_{1,\nu}
    \Phi^0_{\nu}(\vartheta) + \sum_{\nu=1}^{\mathcal{N}^1} c^{1, i+1}_{2,\nu}
    \Phi^1_{\nu}(\vartheta) \quad \text{and} \quad \tilde{\rm u}^{i+1}_{s,\Gamma} \approx \mathfrak{u}^{i+1}_{s,\Gamma},
\end{equation}%
where $\Phi^0_{\nu} = (\phi_{\nu},0)$ and $\Phi^1_{\nu} = (0,
\phi_{\nu})$ are the set of basis functions used to represent the Cauchy
data on the interface $\Gamma$. 

\begin{remark}
{It should be also possible to relax regularity assumption on the domain under consideration. For example, one can consider piecewise-smooth, locally-supported basis functions (defined on $\Gamma$) as the part of the Extension Operator.  For example, \cite{MaguraEtAl2017} use this approach to design a high-order accurate numerical method for the Helmholtz equation, in a geometry with a reentrant corner.  Furthermore, \cite{WoodwardThesis2015,Woodward_2015} combine the DPM together with the XFEM, and design a DPM for linear elasticity in a non-Lipschitz domain (with a cut).}
\end{remark}

{Next, in Appendix~\ref{appendix:DPM:CompleteSystemOfEquations}, we derive} a linear system for the coefficients $(c_{1, \nu}^{s, i+1})_{\nu=1}^{\mathcal{N}^0}$ and $(c_{2, \nu}^{s, i+1})_{\nu=1}^{\mathcal{N}^1}$, by 
{combining the interface conditions (\ref{eqn:PDECD6},~\ref{eqn:PDECD7}), the BEP \eqref{eqn:IntroDPM:BEP}, the Extension Operator (\ref{eqn:ExtOperatorStep1}--\ref{eqn:ExtOperatorStep3}), and the spectral discretization \eqref{eqn:DPM:DiscretizationOfCauchyData}.} 
%
Then, the numerical approximation $u_s^{i+1} \approx u_s(x_j, y_k, t^{i + 1})$ of (\ref{eqn:PDECD1}--\ref{eqn:PDECD7}) at all grid-points $(x_j, y_k) \in N_s^+$ follows directly from the discrete, generalized Green's formula, which we state now. %
\begin{definition}[Discrete, generalized Green's formula.]
    At each time step $t^{i + 1}$, the numerical approximation $u_s^{i+1} \approx u_s(x_j, y_k, t^{i + 1})\rvert_{(x_j, y_k) \in N_s^+}$ of (\ref{eqn:PDECD1}--\ref{eqn:PDECD7}) is given by %
    \begin{equation} \label{eqn:DPM:DGGF}
        u_s^{i + 1} := P_{N_s^+}^{i + 1} u_{\gamma_s}^{i + 1} + G_{\Delta t, h}^s F_s^{i + 1}. 
    \end{equation}%
    Here, $u_{\gamma_s}^{i + 1} = \operatorname{Ex}_s \tilde{\mathfrak{u}}_{s, \Gamma}^{i + 1}$, and $\tilde{\mathfrak{u}}_{s, \Gamma}^{i + 1}$ is constructed from (i) the solution of 
    {the BEP (see Appendix~\ref{appendix:DPM:CompleteSystemOfEquations}) and (ii) the spectral} 
    discretization \eqref{eqn:DPM:DiscretizationOfCauchyData}.  (Recall that $P_{N_s^+}^{i + 1} u_{\gamma_s}^{i + 1}$ 
    {is the Difference Potential of the density $u_{\gamma_s}^{i + 1}$}, 
    while $G_{\Delta t, h}^s F_s^{i + 1}$ 
    {is the Particular Solution}.)
\end{definition}



%
In this work, we also propose a novel feature of DPM, extending the method originally developed in \cite{ryaben2006algorithm} and \cite{phdthesisJA,albright2017high,AlbrightEpshteynSteffen2015,albright2015high,Epshteyn2014,epshteyn2015high} to the composite domain problem (\ref{eqn:PDECD1}--\ref{eqn:PDECD7}) with implicitly-defined geometry. 
The primary difference between 
{Difference Potentials Methods} 
on explicitly-defined \textit{versus} implicitly-defined composite domains is in the approximation of the interface $\Gamma$, which must be done accurately and efficiently, in order to maintain the desired second- or fourth-order accuracy. 

The main idea of 
{DPM-based methods} 
for implicitly-defined geometry is to seek an accurate and efficient explicit parameterization of the implicit boundary/interface. First, we represent the geometry implicitly \textit{via} a level set function $F(x,y)$ on $M^0$. Then we construct a local interpolant $\tilde{F}(x, y)$ of $F(x,y)$ on a subset of $M^0$
near the continuous interface $\Gamma$. Next, we parameterize $\Gamma$ by arc-length using numerical quadrature. With this parameterization, we (i) compute the Fourier series expansion from initial conditions for the Cauchy data 
$\mathfrak{u}_{s, \Gamma}^{i + 1}$
on the implicit interface $\Gamma$, and (ii) construct the extension operators (Definition~\ref{def:DPM:ExtOp}) with $p=2$ or $p=4$. 

\begin{conj}[High-order accuracy of the DPM with implicit geometry] \label{rmk:IntroDPM:Accuracy}
    Due to the second\hyp{} or fourth-order accuracy (in both space and time) of the underlying discretization \eqref{eqn:IntroDPM:FullyDiscreteCD}, the extension operator \eqref{eqn:DPM:ExtOp:DiscreteVersion} with $p=2$ or $p=4$, and the established error estimates and convergence results for the DPM for general linear elliptic boundary value problems 
    on smooth domains (presented in 
    \cite{Reznik1982,Reznik1983Thesis,Ryabenkii2002} and \cite{Godunov2015}), we expect second- and fourth-order accuracy in the maximum norm for the error in the computed solution (\ref{eqn:Measures:ErrSD} or~\ref{eqn:Measures:ErrCD}) for both the single and composite domain parabolic problems.  
\end{conj}

\begin{remark}
    Indeed, in the numerical results (Section~\ref{sec:NumericalResults}) we see that the computed solution \eqref{eqn:DPM:DGGF} at every time level $t^{i+1}$ has accuracy $\mathcal{O}(h^2 + \Delta t^2)$ for the second-order method, and $\mathcal{O}(h^4 + \Delta t^4)$ for the fourth-order method, for both the single and composite domain problems, with explicit or implicit geometry.  See \cite{phdthesisJA,albright2015high,epshteyn2015solution,ryaben2006algorithm} for more details and numerical tests involving explicit (circular and elliptical) geometries. 
\end{remark}

\paragraph{Main Steps of the algorithm:} \label{sec:DPM:MainSteps}
%
Let us summarize the main steps for the Difference Potentials Method. %
\begin{itemize}
    \item \textit{Step 1}: Introduce a computationally simple Auxiliary Domain $\Omega_s^0$ ($s = 1,2$) and formulate the Auxiliary Problem (AP; Definition~\ref{def:IntroDPM:AP}). 
    
    \item \textit{Step 2}: At each time step $t^{i+1}$, compute the Particular Solution $u_{s}^{i+1} = G^{i+1}_{\Delta t, h} F^{i+1}_{s}$, $(x_j,y_k)\in N_s^+$, using the AP with the right-hand side \eqref{eqn:DPM:rhs_PS}. 

    \item \textit{Step 3}: Construct the matrix in the boundary equations \eqref{eqn:DPM:FullSystemUnreduced} (discussed in Appendix~\ref{appendix:DPM:CompleteSystemOfEquations}), derived from the Boundary Equation with Projection (BEP) \eqref{eqn:IntroDPM:BEP}, \textit{via} several solutions of the AP. (When the diffusion coefficients $\lambda_s$ are constant, this is done once, as a pre-processing step before the first time step.) 
    
    \item \textit{Step 4}: Compute the approximation of the density $u^{i+1}_{\gamma_s}$, by applying the Extension Operator \eqref{eqn:DPM:ExtOp:DiscreteVersion} to the solution of \eqref{eqn:DPM:FullSystemUnreduced}. 
    
    \item \textit{Step 5}: Construct the Difference Potentials $P_{N_s^+\gamma_s} u^{i+1}_{\gamma_s}$ of the density $u_{\gamma_s}^{i + 1}$, using the AP with the right-hand side \eqref{eqn:DPM:rhs_DP}. 

    \item \textit{Step 6}: Compute the numerical approximation $u_s^{i+1} \approx u_s(x_j, y_k, t^{i + 1})$ of the PDE (\ref{eqn:PDECD1}--\ref{eqn:PDECD7}) using the discrete, generalized Green's formula \eqref{eqn:DPM:DGGF}. 

\end{itemize}

  \subsection{SBP--SAT--FD} \label{sec:Overview:SBPSAT}


We continue in this section with a brief presentation of SBP--SAT--FD, for solving the parabolic problems presented in Section \ref{sec:Intro:ContinuousProblem}. For more detailed discussions of the SBP--SAT--FD method, we refer the reader to two review papers \cite{Del2014Review,Svard2014}. 

\begin{figure}[tb]
    \centering
    \begin{subfigure}[b]{0.5\linewidth}
        \centering
        \includegraphics[width=\textwidth]{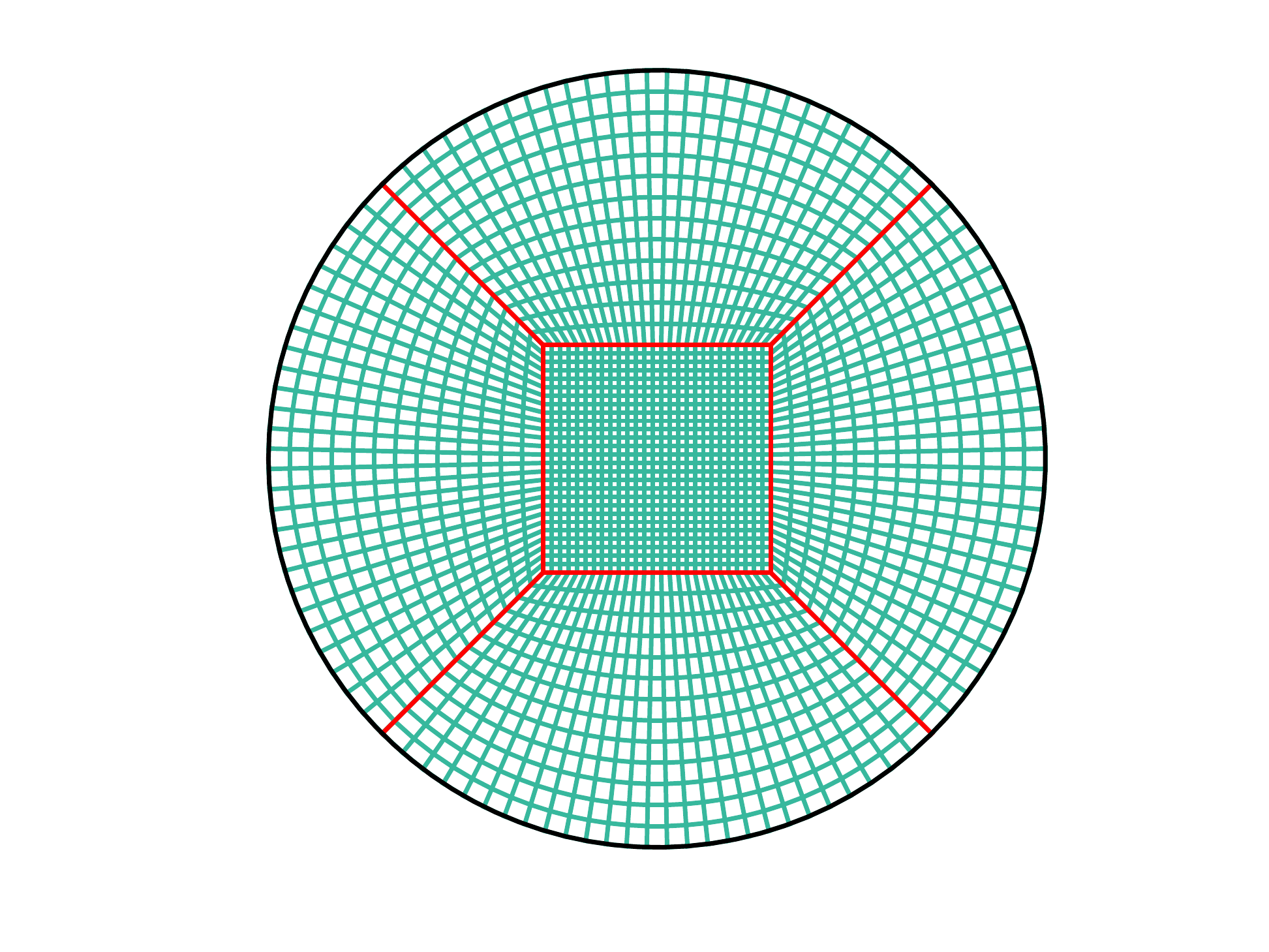}
        \caption{\label{1_Domain_Grid_SBP}}
    \end{subfigure}%
    \begin{subfigure}[b]{0.5\linewidth}
        \centering
        \includegraphics[width=\textwidth]{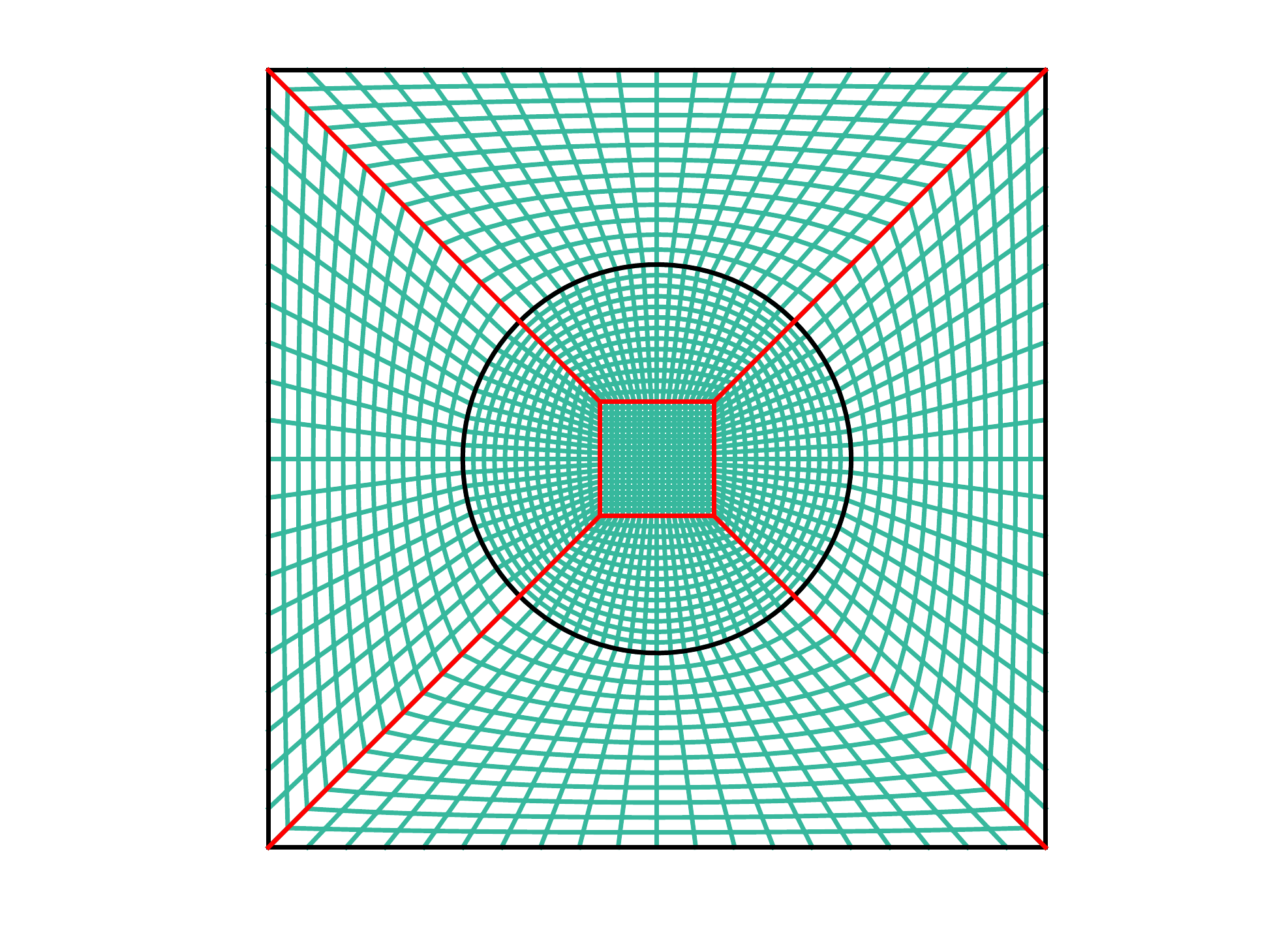}
        \caption{\label{2_Domain_Grid_SBP}}
    \end{subfigure}
    \caption{A \textbf{(a)} circular domain (\textit{cf.} Figure~\ref{fig:CompDomain}), divided into five subdomains; and \textbf{(b)} composite domain, divided into nine subdomains.}
    \label{FD_Grid}
\end{figure}

\siyang{In SBP--SAT--FD, to resolve geometrical features, the physical domain is mapped to a reference domain with a simple geometry, \textit{e.g.}, the unit square.}{The SBP--SAT--FD method was originally used on Cartesian grids. To resolve complex geometries, we consider a grid mapping approach by transfinite interpolation {\cite{Knupp1993}}. }{Comment 10 of reviewer 1} A smooth mapping requires that the physical domain is a quadrilateral, possibly with smooth, curved sides. If the physical domain does not have the desired shape, we then partition the physical domain into subdomains, so that each subdomain can be mapped smoothly to the reference domain. As an example, the single domain of equation ({\ref{eqn:PDESD1}},~{\ref{eqn:PDESD2-3}}), shown in Figure~{\ref{1_Domain_Grid_SBP}}, is divided into five subdomains. The five subdomains consist of one square subdomain, and four identical quadrilateral subdomains (modulo rotation by $\pi/2$) with curved sides. Similarly, the composite domain of equation  ({\ref{eqn:PDECD1}}--{\ref{eqn:PDECD7}}) is divided into nine subdomains, as shown in Figure~{\ref{2_Domain_Grid_SBP}}.  Suitable interface conditions are imposed to patch the subdomains together.

Although the side-length of the centered square is arbitrary (as long as the square is strictly inside the circle), its size and position have a significant impact on the quality of the curvilinear grid. In a high-quality mesh, the elements should not be skewed too much, and the sizes of the elements should be nearly uniform. In practice, it is usually difficult to know \textit{a priori} the optimal way of domain division. 

A Cartesian grid in the reference domain is mapped to a curvilinear grid in each subdomain. The grids are aligned with boundaries and interfaces, thus avoiding small--cut difficulties sometimes associated with embedded methods. \siyang{}{In this paper, we only consider conforming grid interfaces, \textit{i.e.}, the grid points from two adjacent blocks match on the interface. For numerical treatment of non-conforming grid interfaces in the SBP--SAT--FD framework, see {\cite{Kozdon2016,Mattsson2010}.} }{Comment 12 of reviewer 1}




When a physical domain is mapped to a reference domain, the governing equation is transformed to the Cartesian coordinate in the reference domain. The transformed equation is usually in a more complicated form than the original equation. In general, a parabolic problem 
\begin{equation}
u_t = u_{xx}+u_{yy},\quad (x,y)\in\Omega
\end{equation}
in a physical domain will be transformed to 
\begin{equation}\label{heat_vc}
Ju_t = (\alpha u_\xi)_\xi+(\beta u_\eta)_\xi+(\beta u_\xi)_\eta+(\gamma u_\eta)_\eta,\quad  (\xi,\eta)\in [0,1]^2,
\end{equation}
where $(\xi,\eta)$ is the Cartesian coordinate in the unit square, and $J(\xi,\eta)$, $\alpha(\xi,\eta)$, $\beta(\xi,\eta)$, $\gamma(\xi,\eta)$ depend on the geometry of the physical domain and on the chosen mapping. \siyang{}{In particular, we use transfinite interpolation for the grid mapping. In this case, the precise form of \eqref{heat_vc} and the derivation of the grid transformation are presented in Section 3.2 of {\cite{Almquist2014}}.}{Comment 3 from reviewer 2} Even though the original equation is in the simplest form with unit coefficients, the transformed equation has variable coefficients and mixed derivatives. Therefore, it is important to construct multi-block finite difference methods  solving the transformed equation \eqref{heat_vc}. Hence, we need two SBP operators, $D_1\approx\partial/\partial x$ to approximate a first derivative, and $D^{(b)}_2\approx\partial/\partial x(b(x)\partial/\partial x)$ to approximate a second derivative with variable coefficient, where $b(x)>0$ is a known function. Below we discuss SBP properties, and start with the first derivative. 

Consider two smooth functions $u(x), v(x)$ on $x\in [0,1]$. We discretize $[0,1]$ uniformly by $N$ grid points, and denote the restriction of $u(x), v(x)$ onto the grid by $\mathbf{u},\mathbf{v}$, respectively.  Integration by parts states:
\begin{equation}
\int_0^{1} u_xv\ dx= uv\Big|_0^1 - \int_0^{1} uv_x\ dx.
\end{equation}
The SBP operator $D_1$ mimics integration by parts:
 \begin{equation}\label{HD}
 (D_1 \mathbf{u})^T H \mathbf{v}= \mathbf{u}^T B\mathbf{v} - \mathbf{u}^T H D_1\mathbf{v},
 \end{equation}
 where $H$ is symmetric positive definite -- thus defining an inner product -- and 
 $$B=\text{diag}(-1,0,\cdots,0,1).$$ In fact, $H$ is also a quadrature \cite{Del2014}. It is easy to verify that \eqref{HD} is equivalent to 
  \begin{equation}\label{sbp_D1}
D_1^TH+HD_1=B,
 \end{equation}
which is the SBP property for the first derivative operator. At the grid points in the interior of the domain, standard, central, finite-difference stencils can be used in $D_1$, and the weights of the standard, discrete $L_2$-norm are used in $H$. At a few points close to boundaries, special stencils and weights must be constructed in $D_1$ and $H$, respectively, to satisfy \eqref{sbp_D1}. 

The SBP operators $D_1$ were first constructed in \cite{Kreiss1974} and later revisited in \cite{Strand1994}. The SBP norm $H$ can be diagonal or non-diagonal. While non-diagonal norm SBP operators have a better accuracy property than diagonal norm SBP operators, when terms with variable coefficients are present in the equation, a stability proof is only possible with diagonal norm SBP operators. Therefore, we use diagonal norm SBP operators in this paper. 

For a second derivative with variable coefficients, the SBP operators $D^{(b)}_2$ were constructed in \cite{Mattsson2012}. We remark that applying $D_1$ twice also approximates a second derivative, but is less accurate and more computationally expensive than $D^{(b)}_2$.

Due to the choice of centered difference stencils at interior grid points, the order of accuracy of the SBP operators is even at these points, and is often denoted by $2p$. To fulfill the SBP property, at a few grid points near boundaries, the order of accuracy is reduced to $p$ for diagonal norm operators. This detail notwithstanding, such a scheme is often referred to as $2p^{th}$-order accurate. In fact, for the second- and fourth-order SBP--SAT--FD schemes used in this paper to solve parabolic problems, we can expect a second- and fourth-order overall convergence rate, respectively \cite{WangKreiss2016}.

An SBP operator only approximates a derivative. When imposing boundary and interface conditions, it is important that the SBP property is preserved and an energy estimate is obtained. For this reason, we consider the SAT method \cite{Carpenter1994}, where penalty terms are added to the semi-discretization, imposing the boundary and interface conditions weakly. This bears similarities with the Nitsche finite element method \cite{Nitsche1971} and the discontinuous Galerkin method \cite{Hesthaven2008}. 

We note that in \cite{Virta2014}, SBP--SAT--FD methods were developed for the wave equation 
 \begin{equation}\label{wave_vc}
Jv_{tt} = (av_\xi)_\xi+(bv_\eta)_\xi+(bv_\xi)_\eta+(cv_\eta)_\eta,\quad (\xi,\eta)\in [0,1]^2,
\end{equation}
with Dirichlet boundary conditions, Neumann boundary conditions, and interface conditions.  Comparing equation \eqref{wave_vc} with \eqref{heat_vc}, the only difference is that the wave equation has a second derivative in time, while the heat equation has a first derivative in time. The spatial derivatives of \eqref{wave_vc} and \eqref{heat_vc} are the same. 
 
 Assuming homogeneous boundary data for simplified notation, we write the SBP--SAT--FD discretization of \eqref{wave_vc} as 
 \begin{equation}\label{Qv}
 \mathbf{v}_{tt} = Q\mathbf{v},
  \end{equation}
 where $Q$ is the spatial discretization operator including the boundary implementation. For the scheme developed in \cite{Virta2014}, stability is proved by the energy method by multiplying \eqref{Qv} by $\mathbf{v}_t^TH_2$ from the left, 
  \begin{equation}\label{Qv1}
 \mathbf{v}^T_tH_2\mathbf{v}_{tt} = \mathbf{v}^T_tH_2Q\mathbf{v},
  \end{equation}
  where $H_2$ is a diagonal, positive-definite operator, obtained through a tensor product from the corresponding SBP norm, $H$, in one spatial dimension. It is shown in  \cite{Virta2014} that $H_2Q$ is symmetric and negative semi-definite. Therefore, we can write \eqref{Qv1} as 
    \begin{equation*}
\frac{d}{d t} (\mathbf{v}^T_tH_2\mathbf{v}_t - \mathbf{v}^TH_2Q\mathbf{v})=0,
 \end{equation*}
 where the discrete energy, $\mathbf{v}^T_tH_2\mathbf{v}_t - \mathbf{v}^TH_2Q\mathbf{v}$,  for \eqref{wave_vc} is conserved.

If we use the same operator $Q$ to discretize the heat equation \eqref{heat_vc} with the same boundary condition as the wave equation \eqref{wave_vc}, then the scheme
 \begin{equation}\label{Qvheat}
 \mathbf{v}_t = Q\mathbf{v},
  \end{equation}
 is also stable. To see this, we multiply \eqref{Qvheat} by $\mathbf{v}^TH_2$ from the left, and obtain
\begin{equation}\label{energy_wave}
\frac{d}{d t} (\mathbf{v}^TH_2\mathbf{v}) = \mathbf{v}^TH_2Q\mathbf{v} \leq 0,
\end{equation}
 where $\mathbf{v}^TH_2\mathbf{v}$ is the discrete energy for \eqref{heat_vc}. In this paper, we use the spatial discretization operators developed in \cite{Virta2014} to solve both the single (\ref{eqn:PDESD1},~\ref{eqn:PDESD2-3}) and composite domain problems (\ref{eqn:PDECD1}--\ref{eqn:PDECD7}).
 
In \cite{Berg2012}, SBP--SAT--FD methods are discussed for the one-dimensional heat equation with constant coefficients, both in a single domain and a composite domain. In theory, these schemes can also be generalized to solve equation \eqref{heat_vc}, but are different from the ones used in this paper.

  %
  %
  %
  %
  %
  %

\section{Test Problems} \label{sect:TestProblems}

In this section, we first list the test problems that we will consider (in Section~\ref{sec:TestProblems:List}), and then briefly motivate and discuss these choices (in Section~\ref{sec:TestProblems:Motivation}).  The tests we propose are ``manufactured solutions'', in the sense that we state an exact solution $u(x, y, t)$ or $(u_1(x, y, t), \, u_2(x, y, t))$ and a diffusion coefficient $\lambda(t)$ or $(\lambda_1, \lambda_2)$.  From (\ref{eqn:PDESD1},~\ref{eqn:PDESD2-3}) (for the single domain problem) or (\ref{eqn:PDECD1}--\ref{eqn:PDECD7}) (for the composite domain problem) we compute the (i) right-hand side, (ii) initial conditions, (iii) boundary condition, and (iv) functions $(\mu_1(x, y, t), \, \mu_2(x, y, t))$ for the interface/matching conditions.  Then, (i--iv), together with the diffusion coefficient, serve as the inputs for our numerical methods. 

\subsection{List of test problems} \label{sec:TestProblems:List}

\begin{enumerate}

  \item Single-domain, with an explicitly-defined boundary for DPM and SBP--SAT--FD, or an implicitly-defined boundary for cut--FEM. 

  \begin{enumerate}
    \item Constant diffusion (Test Problem 1A; TP--1A):  Consider the PDE (\ref{eqn:PDESD1},~\ref{eqn:PDESD2-3}), with $\lambda(t) \equiv 1$, $\Omega = \{ (x, y) \in \RR^2 : x^2 + y^2 \leq 1 \}$, and the final time $T=1.0$.  Then, TP--1A (adapted from \cite{albright2015high}), is given by %
    \begin{equation} \label{eqn:TP1A}
      u(x, y, t) = x^9 y^8 e^{-t}. \tag{TP--1A}
    \end{equation}%

    \item Time-varying diffusion (Test Problem 3A; TP--3A):  Same as TP--1A, but with diffusion coefficient 
    \begin{equation}\label{eqn:TP3A}
    \lambda(t) = 11/10 + \sin (10 \pi t). \tag{TP--3A}
    \end{equation} 

  \end{enumerate}

  \item Composite-domain, with an explicitly-defined interface (for DPM and SBP--SAT--FD) or implicitly-defined interface (for cut--FEM and DPM).  Consider the PDE (\ref{eqn:PDECD1}--\ref{eqn:PDECD7}), with $\Omega = [-2, 2] \times [-2, 2]$, $\Omega_2 = \{ (x, y) \in \RR^2 : x^2 + y^2 \leq 1 \}$, $\Omega_1 = \Omega \setminus \Omega_2$, $\Gamma = \{ (x, y) \in \RR^2 : x^2 + y^2 = 1 \}$, and the final time $T = 1.0$. %
  \begin{enumerate}
    \item (Test Problem 2A; TP--2A):  A modified version of the test adapted from \cite{albright2015high,li2006immersed}.  Let $(\lambda_1, \lambda_2) = (10, 1)$, and %
    \begin{equation} \label{eqn:TP2A}
      u(x, y, t) = \begin{cases}
        e^{-t} \sin x \cos y, & (x, y) \in \Omega_1, \\
        e^{-t} (x^2 - y^2), & (x, y) \in \Omega_2. 
      \end{cases} \tag{TP--2A}
    \end{equation}

    \item High-frequency oscillations (Test Problem 2B; TP--2B):  A modified version of the test adapted from \cite{albright2015high}.  Let $(\lambda_1, \lambda_2) = (10, 1)$, and %
    \begin{equation} \label{eqn:TP2B}
      u(x, y, t) = \begin{cases}
        e^{-t} \sin (3 \pi x) \cos (7 \pi y), & (x, y) \in \Omega_1, \\
        e^{-t} (x^2 - y^2), & (x, y) \in \Omega_2. 
      \end{cases} \tag{TP--2B}
    \end{equation}

    \item Large contrast in diffusion coefficients, and large jumps in both solution and flux at interface (Test Problem 2C; TP--2C):  Let $(\lambda_1, \lambda_2) = (1000, 1)$, and %
    \begin{equation} \label{eqn:TP2C}
      u(x, y, t) = \begin{cases}
        0, & (x, y) \in \Omega_1, \\
        1000 \sin(10 t) x^4 y^5, & (x, y) \in \Omega_2. 
      \end{cases} \tag{TP--2C}
    \end{equation}%

  \end{enumerate}%

\end{enumerate}

\subsection{Motivation of the chosen test problems} \label{sec:TestProblems:Motivation}




Test Problem 1A \eqref{eqn:TP1A} involves a high-degree polynomial, with total degree of $17$.  This is a rather straightforward test problem, which allows us to establish a good ``baseline'' with which to compare each method.  The choice of high degree ensures that there will be no cancellation of local truncation error, so that we should see -- at most -- second- or fourth-order convergence for the given methods, barring some type of superconvergence.
Next, \eqref{eqn:TP3A} adds on (incrementally) the complication of time-varying diffusion. 

Likewise, \eqref{eqn:TP2A} offers a straightforward ``baseline'' with which to consider the interface problem:  The test problem is piecewise-smooth, and the geometry is simplified (see Remark~\ref{rmk:MotivationForGeometries}).  However, there is a jump in both the analytical solution and its flux, which requires a well-designed numerical method to accurately approximate.  
Moreover, \eqref{eqn:TP2A} was first proposed in \cite{li2006immersed} (see also \cite{albright2015high}), and is a good comparison with the immersed interface method therein. 

Then, \eqref{eqn:TP2B} adds additional challenges onto \eqref{eqn:TP2A} in the form of much higher-frequency oscillations; while \eqref{eqn:TP2C} adds onto \eqref{eqn:TP2A} in the form of both (i) large contrast in diffusion, and (ii) large jumps in the analytical solution and its flux.



\section{Numerical results} \label{sec:NumericalResults}

  \subsection{Time discretization} \label{sect:TimeStepping}

The spatial discretization for each method is discussed in Section~\ref{sec:Intro:OverviewMethods}. 
For the time discretization, the backward differentiation formulas of second- and fourth-order (BDF2 and BDF4) are used for the second- and fourth-order methods, respectively. In each case, the time-step is given by %
\begin{equation} \label{eqn:time-step}
    \Delta t=0.5h.
\end{equation}%
However, note that $h$ in \eqref{eqn:time-step} bears different physical meanings for each method. Indeed, for cut--FEM, $h$ is the average distance between the Gauss--Lobatto points; for DPM, $h$ is the grid spacing in the uniform, Cartesian grid $M^0$ (see the text prior to \eqref{eqn:IntroDPM:FullyDiscreteCD}); and for SBP--SAT--FD, $h$ is the minimum grid spacing in the reference domain.

  \subsection{Measure for comparison} \label{sect:MeasuresForComparison}


Let $\smash{u_{j, k}^i}$ denote the computed numerical approximation of $u(x, y, t)$ at the grid-point $(x_j, y_k) \in \Omega$ and time $t^i = i \Delta t \in (0, T]$.
For the three methods, we will compare the size of the maximum error in $u$ at the grid points, with respect to the number of degrees of freedom (DOF).
For the single domain problem (\ref{eqn:PDESD1},~\ref{eqn:PDESD2-3}), the maximum error is computed as:
\begin{align} 
E & := \max_{t^i \in (0, T]} \max_{(x_j, y_k) \in \Omega} \abs{u(x_j, y_k, t^i) - u_{j, k}^i}, \label{eqn:Measures:ErrSD} \\
%
    \intertext{and for the composite domain problem (\ref{eqn:PDECD1}--\ref{eqn:PDECD7}) as:}
%
E & := \max_{t^i \in (0, T]} \max_{(x_j, y_k) \in \Omega_1 \cup \Omega_2} \abs{u(x_j, y_k, t^i) - u_{j, k}^i}. \label{eqn:Measures:ErrCD}
\end{align}%




  \subsection{Convergence results} \label{sect:ConvergenceResults}


\renewcommand{\floatpagefraction}{.9}

In the following tables and figures, we state the number of degrees of freedom in the grid, maximum error (\ref{eqn:Measures:ErrSD}, \ref{eqn:Measures:ErrCD} for the single- and composite-domain problems, respectively), and an estimate of the rate of convergence.


In Tables~\ref{tbl:combined-results-TP1A}--\ref{tbl:combined-results-TP2C}, the estimate of rate of convergence is computed as follows.  Let $(\text{DOF}_{n}, \, E_n)$ be given, with $n = 1, 2, 3$ referring to the first, second, and third grids 
(from coarsest to finest).  Then, for $n = 2, 3$, compute the standard estimate %
\begin{equation} \label{eqn:ConvEst}
    \rho_n = \frac{\log(E_{n - 1} / E_n)}{\log(\text{DOF}_{n - 1} / \text{DOF}_{n})}, 
\end{equation}%
which is the estimated rate of convergence, denoted in Tables~\ref{tbl:combined-results-TP1A}--\ref{tbl:combined-results-TP2C} by ``Rate''.

\begin{figure}[h]
    \centering
    \includegraphics[scale=0.85,page=1]{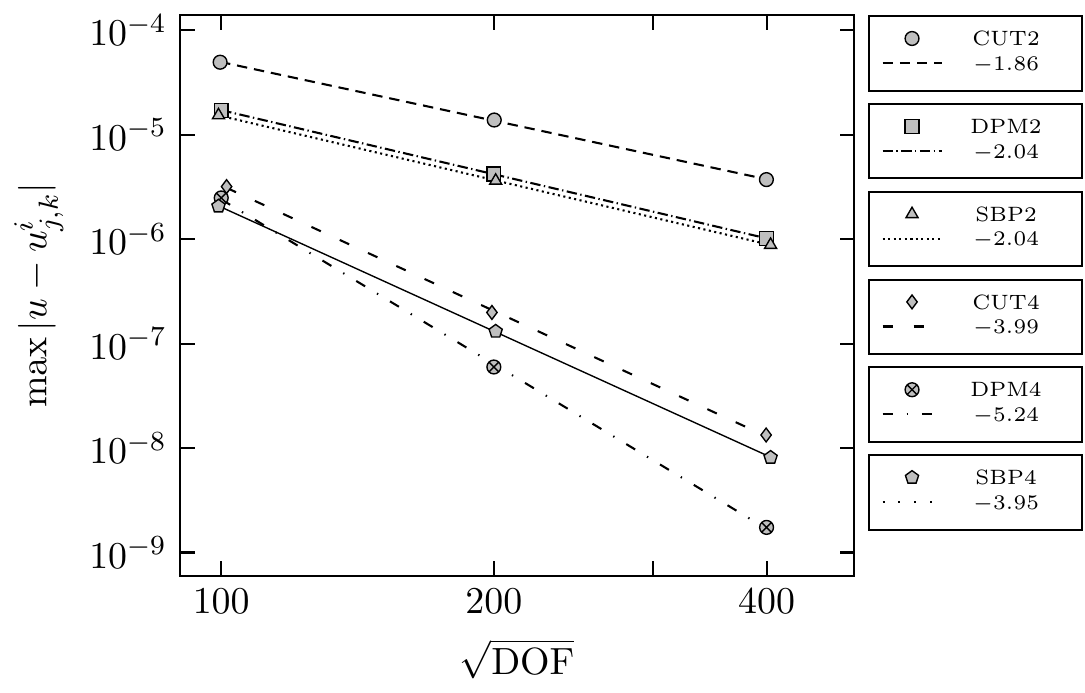}
    \caption{Log--log plot of absolute error \eqref{eqn:Measures:ErrSD} \textit{versus} $\sqrt{\text{DOF}}$, and estimated rate of convergence, for the second- and fourth-order versions of each method, applied to Test Problem 1A \eqref{eqn:TP1A}.  See Table~\ref{tbl:combined-results-TP1A} for more details.}
    \label{fig:combined-results-TP1A}
\end{figure}

\begin{figure}[h]
    \centering
    \includegraphics[scale=0.85,page=1]{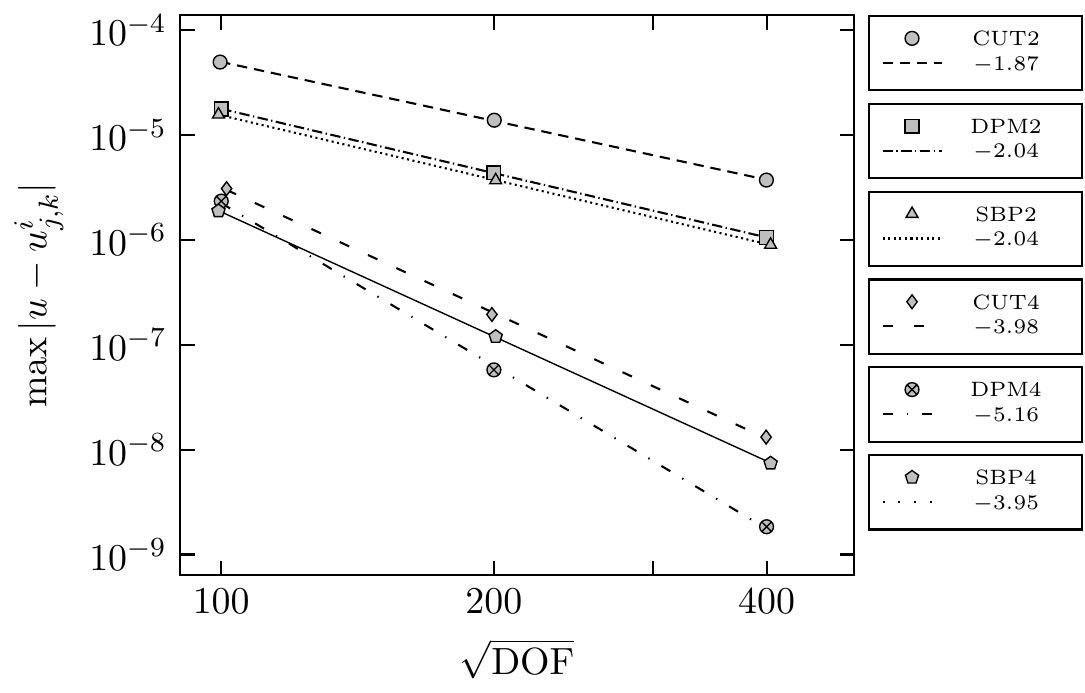}
    \caption{Log--log plot of absolute error \eqref{eqn:Measures:ErrSD} \textit{versus} $\sqrt{\text{DOF}}$, and estimated rate of convergence, for the second- and fourth-order versions of each method, applied to Test Problem 3A \eqref{eqn:TP3A}.  See Table~\ref{tbl:combined-results-TP3A} for more details.}
    \label{fig:combined-results-TP3A}
\end{figure}

\begin{figure}[h]
    \centering
    \includegraphics[scale=0.85,page=1]{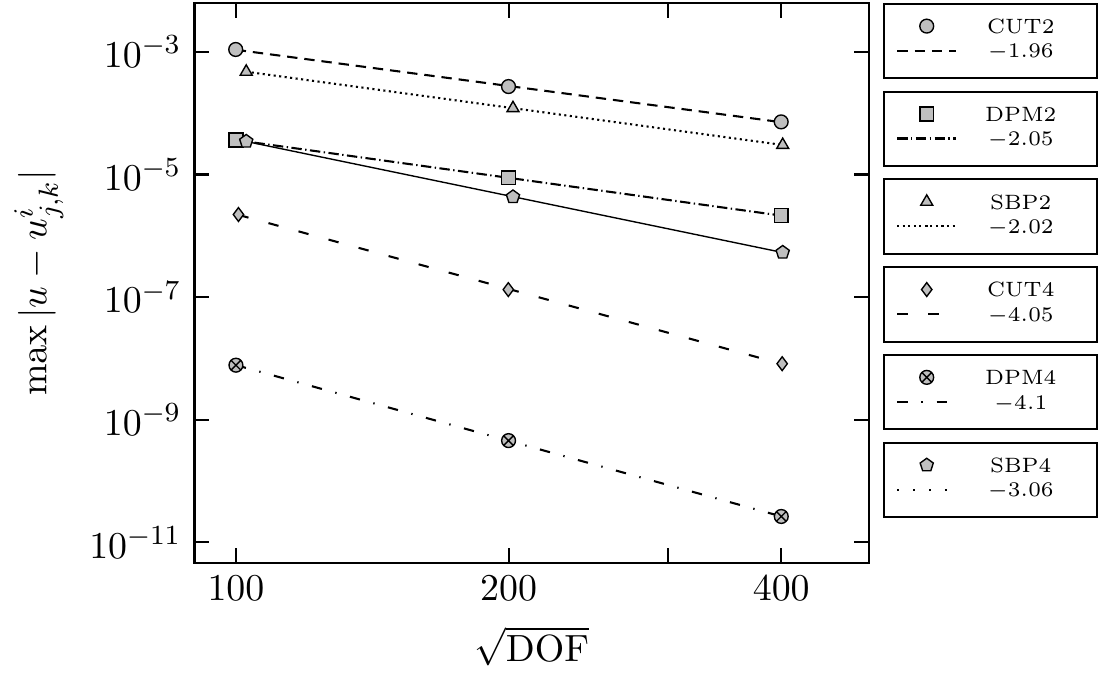}
    \caption{Log--log plot of absolute error \eqref{eqn:Measures:ErrCD} \textit{versus} $\sqrt{\text{DOF}}$, and estimated rate of convergence, for the second- and fourth-order versions of each method (cut--FEM, DPM with explicit geometry, SBP--SAT--FD), applied to Test Problem 2A \eqref{eqn:TP2A}.  See Table~\ref{tbl:combined-results-TP2A} for more details.} 
    \label{fig:combined-results-TP2A}
\end{figure}

\begin{figure}[h]
    \centering
    \includegraphics[scale=0.85,page=1]{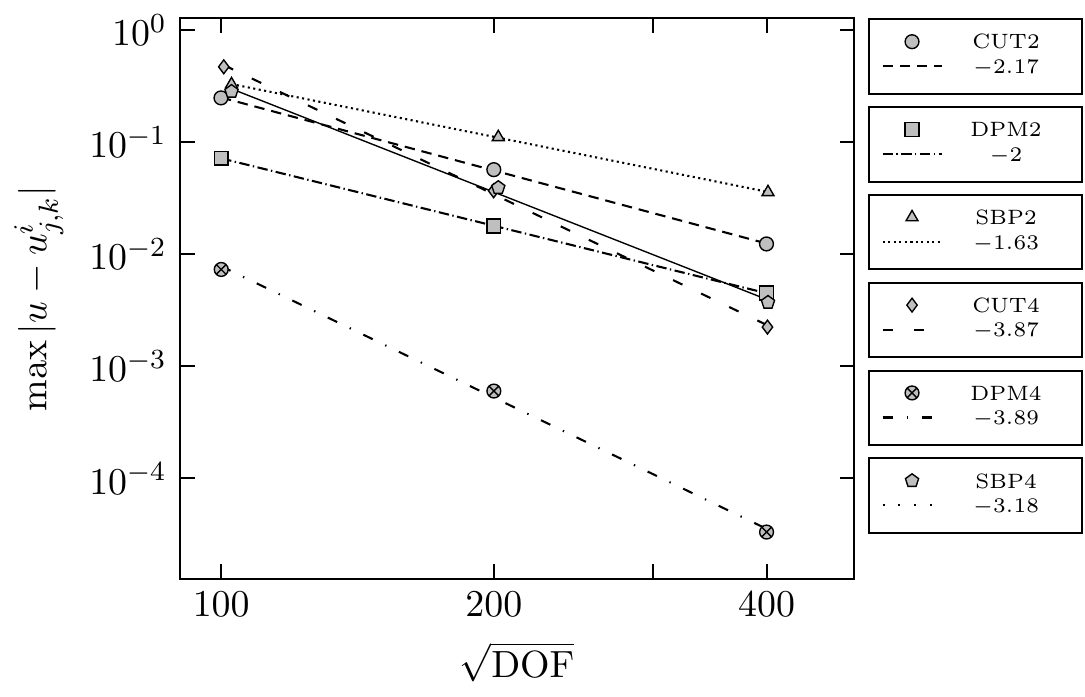}
    \caption{Log--log plot of absolute error \eqref{eqn:Measures:ErrCD} \textit{versus} $\sqrt{\text{DOF}}$, and estimated rate of convergence, for the second- and fourth-order versions of each method (cut--FEM, DPM with explicit geometry, SBP--SAT--FD), applied to Test Problem 2B \eqref{eqn:TP2B}.  See Table~\ref{tbl:combined-results-TP2B} for more details.} 
    \label{fig:combined-results-TP2B}
\end{figure}

\begin{figure}[h]
    \centering
    \includegraphics[scale=0.85,page=1]{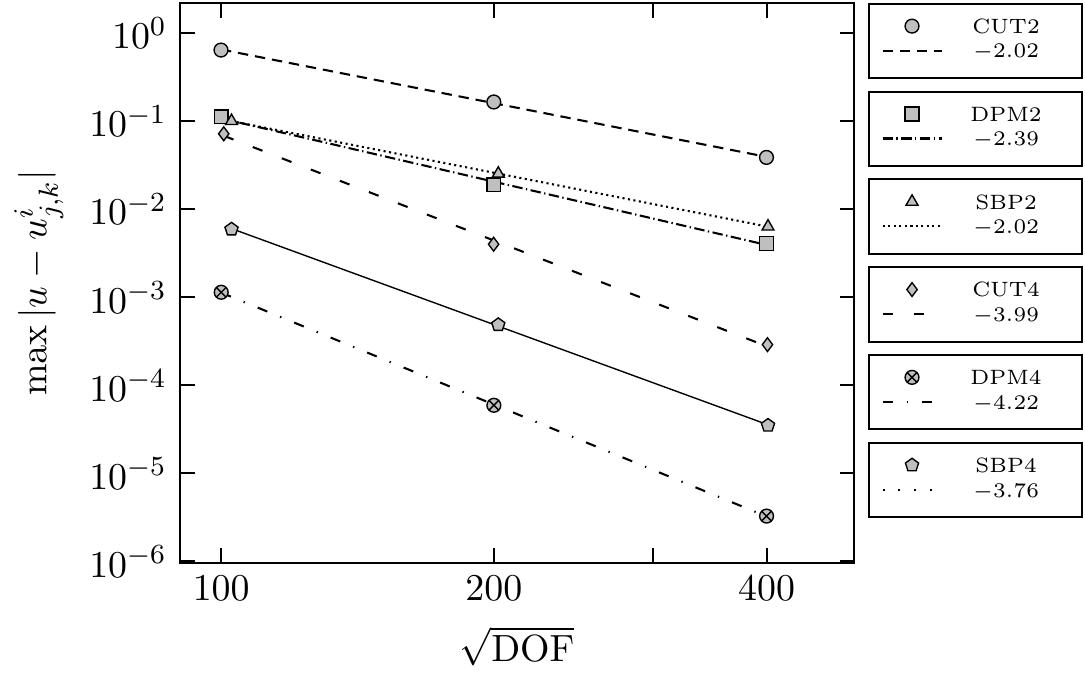}
    \caption{Log--log plot of absolute error \eqref{eqn:Measures:ErrCD} \textit{versus} $\sqrt{\text{DOF}}$, and estimated rate of convergence, for the second- and fourth-order versions of each method (cut--FEM, DPM with explicit geometry, SBP--SAT--FD), applied to Test Problem 2C \eqref{eqn:TP2C}.  See Table~\ref{tbl:combined-results-TP2C} for more details.} 
    \label{fig:combined-results-TP2C}
\end{figure}

In Figures~\ref{fig:combined-results-TP1A},~\ref{fig:combined-results-TP3A},~\ref{fig:combined-results-TP2A}--\ref{fig:combined-results-TP2C}, the estimate of rate of convergence is computed differently.
Computing a least-square linear regression for the data $(\log_{10}(\sqrt{\text{DOF}_n})$, $\log_{10}(E_n))$ 
gives a line with slope $m$, where $m$ is the estimate of rate of convergence, reported in the legend on the right side of each figure.

\begin{table}[h]

\centering 
\footnotesize 
\sisetup{
  output-exponent-marker = \text{ E},
  exponent-product={},
  retain-explicit-plus,
  group-digits=integer,
  group-separator={,},
  group-minimum-digits=4,
  table-number-alignment=center}

\begin{tabular}{S[table-format=6]S[table-format=+1.4e+2]S[table-format=1.2]S[table-format=6]S[table-format=+1.4e+2]S[table-format=1.2]}
\toprule 
{DOF} & {$E$: CUT2} & {Rate}  & {DOF} & {$E$: CUT4} & {Rate} \\ 
\midrule 
9944 & 4.9327e-05 & {---} & 10276 & 3.1799e-06 & {---} \\ 
40072 & 1.3798e-05 & 1.80  & 39613 & 1.9848e-07 & 4.00  \\ 
159912 & 3.7114e-06 & 1.88  & 159700 & 1.3330e-08 & 3.82  \\ 
\midrule[\heavyrulewidth] 
{DOF} & {$E$: DPM2} & {Rate}  & {DOF} & {$E$: DPM4} & {Rate} \\ 
\midrule 
10000 & 1.7105e-05 & {---} & 10000 & 2.4782e-06 & {---} \\ 
40000 & 4.1980e-06 & 2.03  & 40000 & 5.9672e-08 & 5.38  \\ 
160000 & 1.0135e-06 & 2.05  & 160000 & 1.7396e-09 & 5.10  \\ 
\midrule[\heavyrulewidth] 
{DOF} & {$E$: SBP2} & {Rate}  & {DOF} & {$E$: SBP4} & {Rate} \\ 
\midrule 
9861 & 1.5328e-05 & {---} & 9861 & 2.0636e-06 & {---} \\ 
40365 & 3.6210e-06 & 2.08  & 40365 & 1.3083e-07 & 3.98  \\ 
163317 & 8.8008e-07 & 2.04  & 163317 & 8.1180e-09 & 4.01  \\ 
\bottomrule 
\end{tabular}

\caption[
]{
  Convergence in the maximum norm \eqref{eqn:Measures:ErrSD}, for the second- and fourth-order versions of each method, applied to Test Problem 1A \eqref{eqn:TP1A}, 
  with diffusion coefficient $\lambda = 1$, and time-step $\Delta t = 0.5 h$. 
  
}

\label{tbl:combined-results-TP1A}

\end{table}

Overall, we see in Tables~\ref{tbl:combined-results-TP1A}--\ref{tbl:combined-results-TP2C} that the error for second-order methods (denoted, for brevity, as CUT2, DPM2, SBP2) on the finest mesh is similar, or sometimes larger, than the error for fourth-order methods (denoted CUT4, DPM4, SBP4) on the coarsest mesh -- this illustrates the effectiveness of higher-order methods, when high accuracy is important. 
Additionally, comparing the three methods together, the size of the errors for the single-domain problems (\ref{eqn:TP1A}, \ref{eqn:TP3A}) are similar, up to a constant factor; while for the composite-domain problems (\ref{eqn:TP2A}, \ref{eqn:TP2B}, \ref{eqn:TP2C}) we do see differences of one or two orders of magnitude, with the DPM having the smallest errors.

In Table~\ref{tbl:combined-results-TP1A} and Figure~\ref{fig:combined-results-TP1A}, we observe that the measured rates of convergence for the numerical approximation of Test Problem 1A \eqref{eqn:TP1A} are all $\approx 2$ (for the second-order versions) or $\approx 4$ (for the fourth-order versions), except for DPM4, which for this test problem is superconvergent, with fifth-order convergence.  Such higher-than-expected convergence might occur due to several reasons -- for example, (i) if the geometry is smooth; (ii) if the magnitude of the derivatives have fast decay (effectively reducing the local truncation error by a factor of $h$); or (iii) if there is cancellation of error due to symmetries in the geometry, or in the analytical solution.

\begin{table}[h]

\centering 
\footnotesize 
\sisetup{
  output-exponent-marker = \text{ E},
  exponent-product={},
  retain-explicit-plus,
  group-digits=integer,
  group-separator={,},
  group-minimum-digits=4,
  table-number-alignment=center}

\begin{tabular}{S[table-format=6]S[table-format=+1.4e+2]S[table-format=1.2]S[table-format=6]S[table-format=+1.4e+2]S[table-format=1.2]}
\toprule 
{DOF} & {$E$: CUT2} & {Rate}  & {DOF} & {$E$: CUT4} & {Rate} \\ 
\midrule 
9944 & 4.9605e-05 & {---} & 10276 & 3.0791e-06 & {---} \\ 
40072 & 1.3851e-05 & 1.80  & 39613 & 1.9435e-07 & 3.99  \\ 
159912 & 3.7176e-06 & 1.89  & 159700 & 1.3161e-08 & 3.81  \\ 
\midrule[\heavyrulewidth] 
{DOF} & {$E$: DPM2} & {Rate}  & {DOF} & {$E$: DPM4} & {Rate} \\ 
\midrule 
10000 & 1.7721e-05 & {---} & 10000 & 2.3422e-06 & {---} \\ 
40000 & 4.3619e-06 & 2.02  & 40000 & 5.7588e-08 & 5.35  \\ 
160000 & 1.0526e-06 & 2.05  & 160000 & 1.8398e-09 & 4.97  \\ 
\midrule[\heavyrulewidth] 
{DOF} & {$E$: SBP2} & {Rate}  & {DOF} & {$E$: SBP4} & {Rate} \\ 
\midrule 
9861 & 1.5665e-05 & {---} & 9861 & 1.8858e-06 & {---} \\ 
40365 & 3.6965e-06 & 2.08  & 40365 & 1.1949e-07 & 3.98  \\ 
163317 & 8.9731e-07 & 2.04  & 163317 & 7.4149e-09 & 4.01  \\ 
\bottomrule 
\end{tabular}

\caption[
]{
  Convergence in the maximum norm \eqref{eqn:Measures:ErrSD}, for the second- and fourth-order versions of each method, applied to Test Problem 3A \eqref{eqn:TP3A}, 
  with diffusion coefficient $\lambda(t) = 1.1 + \sin(\pi t)$, and time-step $\Delta t = 0.5h$. 
  
}

\label{tbl:combined-results-TP3A}

\end{table}

Table~\ref{tbl:combined-results-TP3A} and Figure~\ref{fig:combined-results-TP3A} show the numerical results for \eqref{eqn:TP3A}.  This test problem has the same manufactured solution as \eqref{eqn:TP1A}, but with a time-varying diffusion coefficient.  Despite this added complexity, the numerical results are the same order of accuracy, and in many cases the errors are the same up to seven digits, when compared with the results for \eqref{eqn:TP1A}.  \redtext{This similarity in the numerical results demonstrates that the three methods can robustly handle time-varying diffusion coefficients.}{Comment 16 Reviewer 1}

\begin{figure}[p]
    \centering
    \begin{subfigure}[b]{0.9\linewidth}
        \centering
        \includegraphics[width=3.5in]{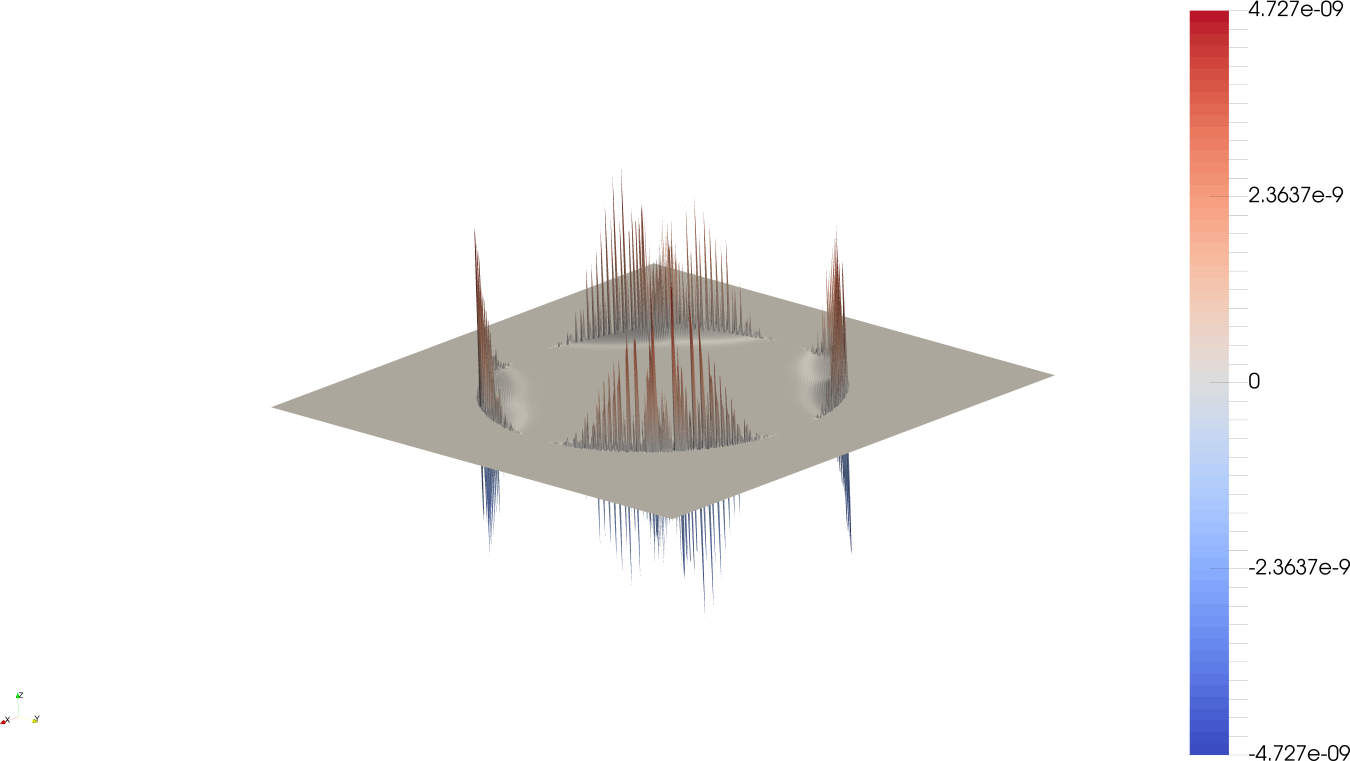}
        \caption{}
        \label{fig:error_CUT_3A}
    \end{subfigure}

    \begin{subfigure}[b]{0.9\linewidth}
        \centering
        \includegraphics[width=3.5in]{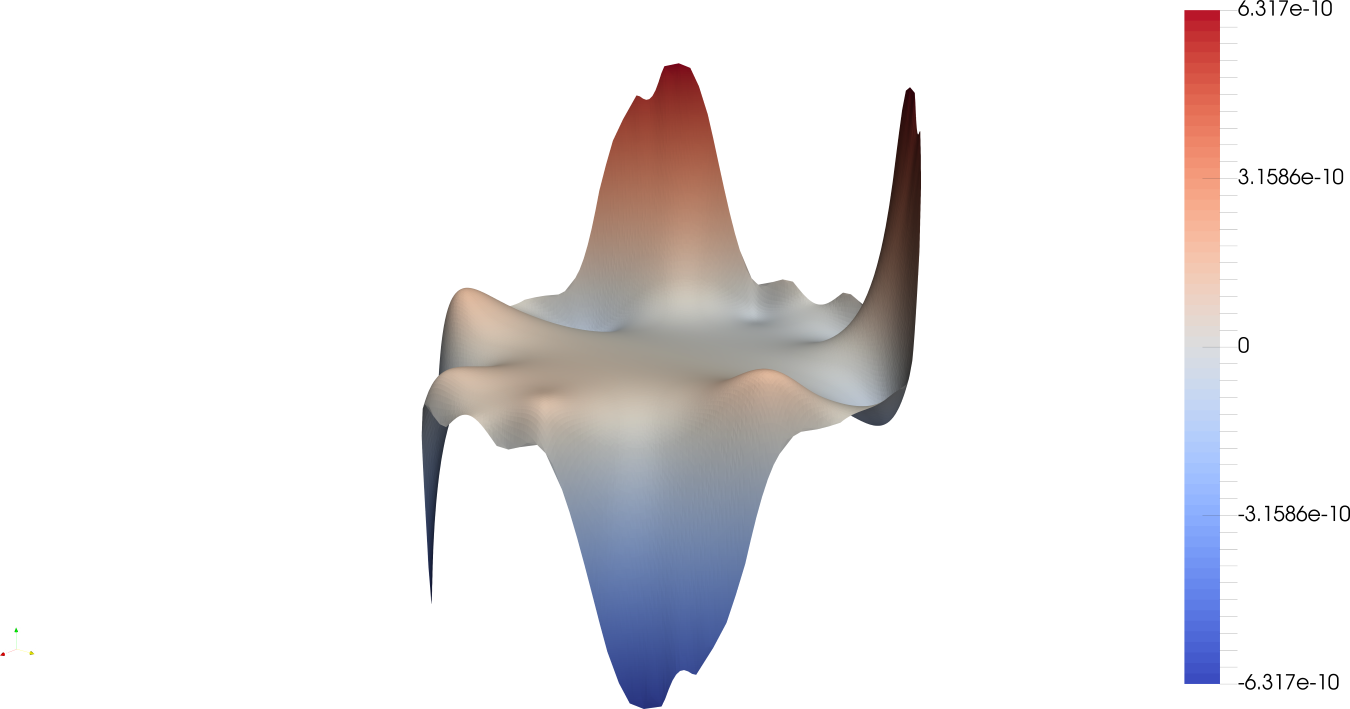}
        \caption{}
        \label{fig:error_DPM_TP3a400}
    \end{subfigure}

    \begin{subfigure}[b]{0.9\linewidth}
        \centering
        \includegraphics[width=3.5in]{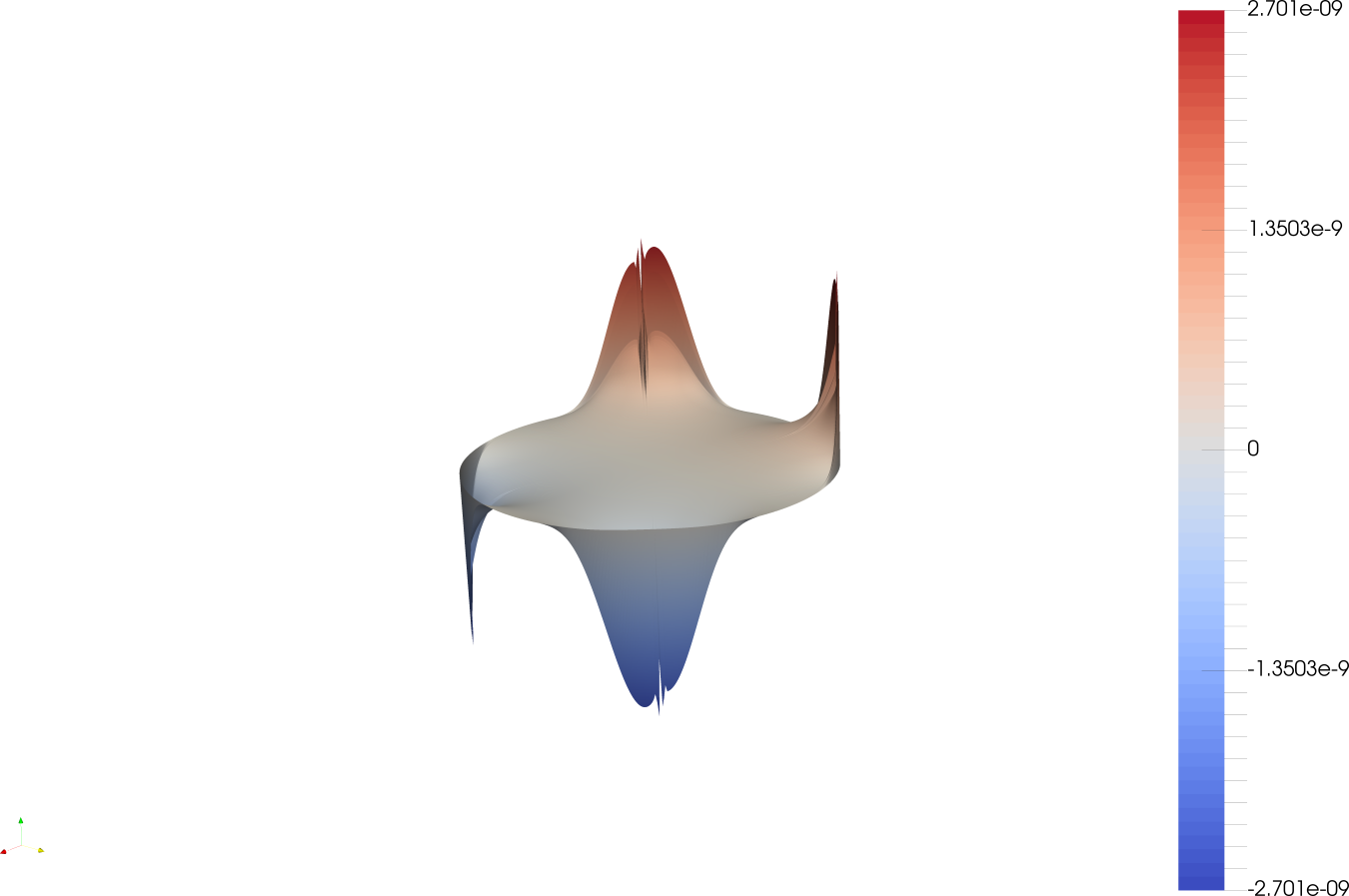}
        \caption{}
        \label{fig:error_FD_3A}
    \end{subfigure}

    \caption{Plot of error at the final time $T = 1.0$, for the fourth-order versions of \textbf{(a)} cut--FEM, \textbf{(b)} DPM, and \textbf{(c)} SBP--SAT--FD, respectively, applied to Test Problem 3A \eqref{eqn:TP3A}.}
    \label{fig:PlotsErrTP3A}
\end{figure}

\begin{figure}[p]
    \centering
    \begin{subfigure}[b]{0.9\linewidth}
        \centering
        \includegraphics[width=3.5in]{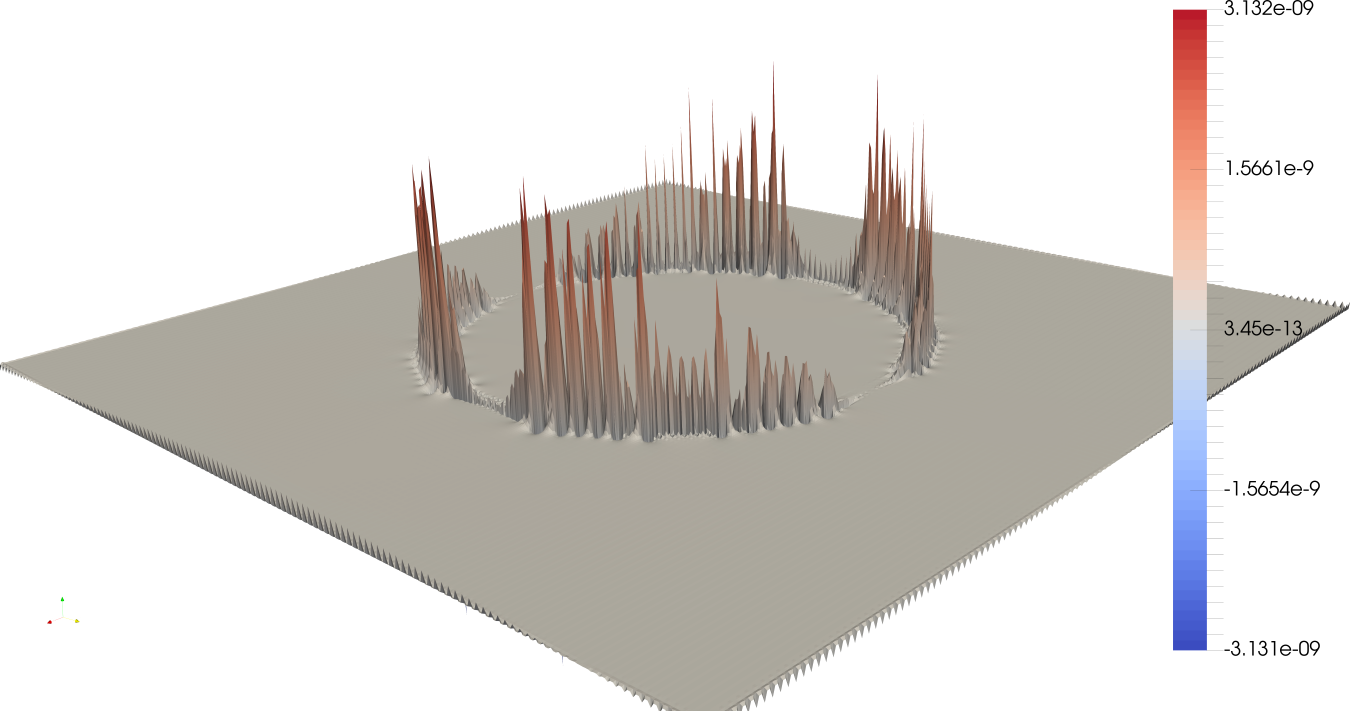}
        \caption{}
        \label{fig:error_CUT_4A_fine}
    \end{subfigure}

    \begin{subfigure}[b]{0.9\linewidth}
        \centering
        \includegraphics[width=3.5in]{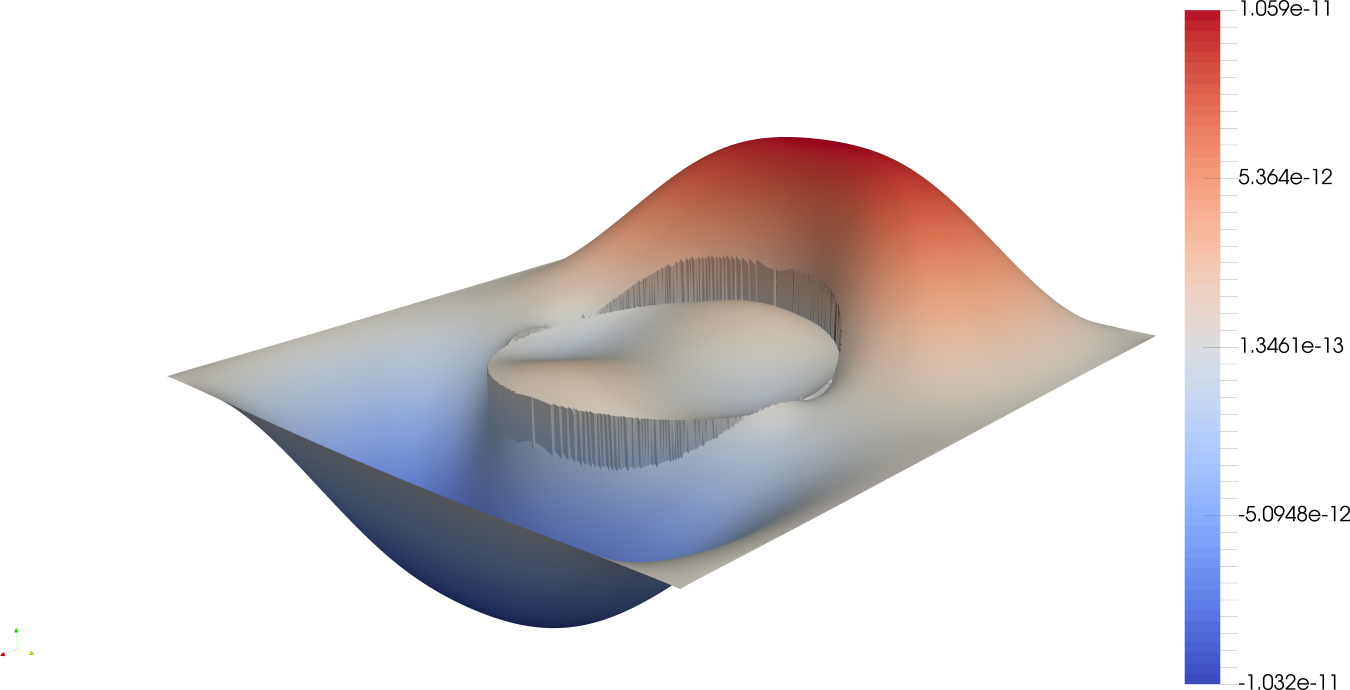}
        \caption{}
        \label{fig:error_DPM_TP2a400}
    \end{subfigure}

    \begin{subfigure}[b]{0.9\linewidth}
        \centering
        \includegraphics[width=3.5in]{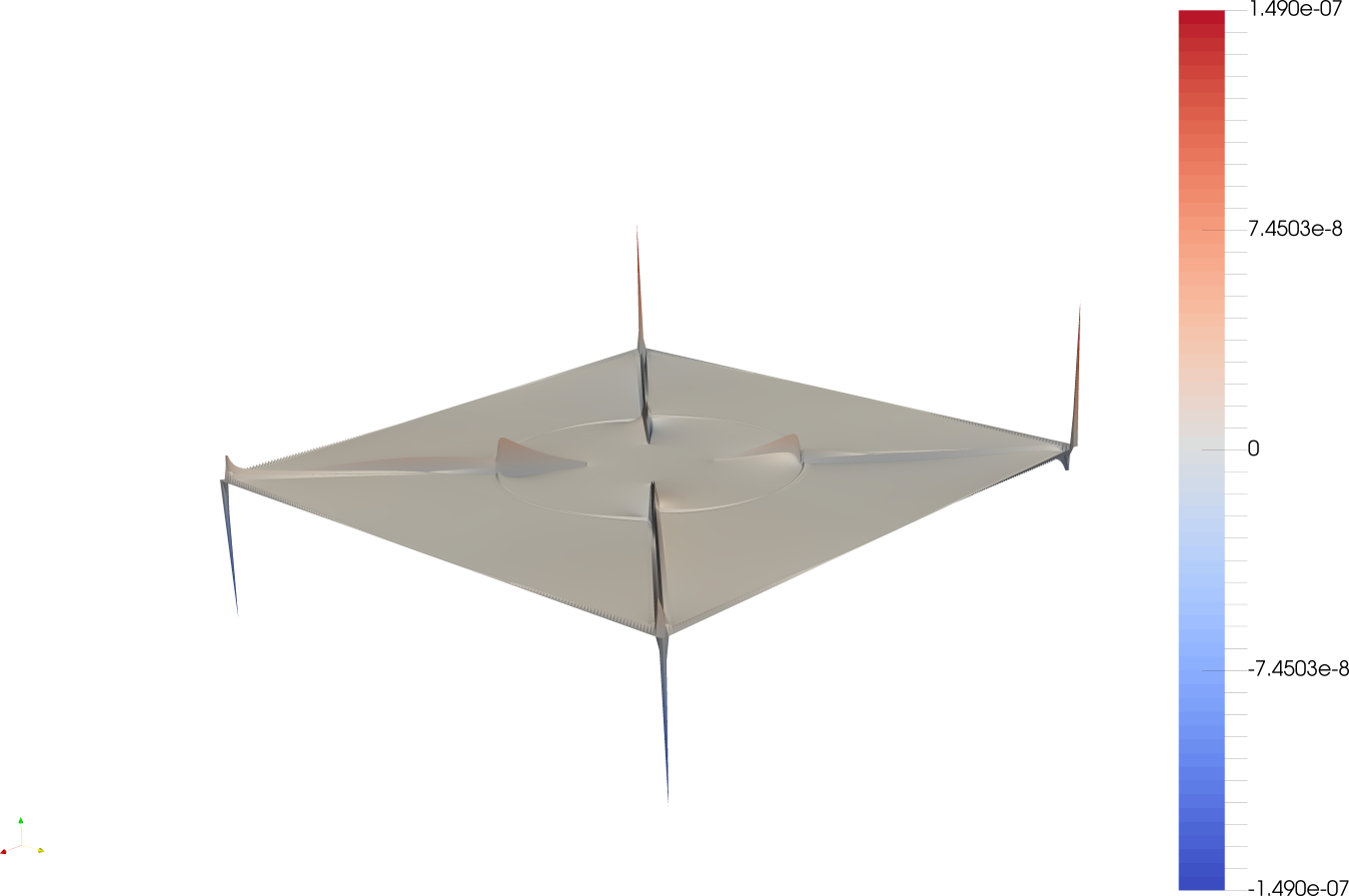}
        \caption{}
        \label{fig:error_FD_2A}
    \end{subfigure}

    \caption{Plot of error at the final time $T = 1.0$, for the fourth-order versions of \textbf{(a)} cut--FEM, \textbf{(b)} DPM, and \textbf{(c)} SBP--SAT--FD, respectively, applied to Test Problem 2A \eqref{eqn:TP2A}.}
    \label{fig:PlotsErrTP2A}
\end{figure}

The plots of spatial error at the final time $T = 1.0$, shown in Figure~\ref{fig:PlotsErrTP3A}, are representative of other tests (not included in this text) on a single circular domain.  The error in the cut--FEM solution presents largely at the boundary; the error in the DPM solution typically has smooth error, even for grid points very near $\Gamma$; while the error in the SBP--SAT--FD solution is not smooth at interfaces introduced by the domain partitioning.


The plots of spatial error at the final time $T =1.0$ for \eqref{eqn:TP2A} are shown in Figure~\ref{fig:PlotsErrTP2A}.  These plots are fairly representative of the other composite domain tests reported herein, and also of others test problems not included in this work.  As in Figure~\ref{fig:PlotsErrTP3A}, the cut--FEM has its largest error at degrees of freedom on cut (intersected) elements; the DPM has piecewise smooth error, including even grid points at the boundary/interface; and the SBP--SAT--FD has its largest error at the interfaces between computational subdomains, with particularly pronounced error at the corners of $\Omega$, where the grid is most stretched.

\begin{table}[h]

\centering 
\footnotesize 
\sisetup{
  output-exponent-marker = \text{ E},
  exponent-product={},
  retain-explicit-plus,
  group-digits=integer,
  group-separator={,},
  group-minimum-digits=4,
  table-number-alignment=center}

\begin{tabular}{S[table-format=6]S[table-format=+1.4e+2]S[table-format=1.2]S[table-format=6]S[table-format=+1.4e+2]S[table-format=1.2]}
\toprule 
{DOF} & {$E$: CUT2} & {Rate}  & {DOF} & {$E$: CUT4} & {Rate} \\ 
\midrule 
9988 & 1.0933e-03 & {---} & 10129 & 2.2215e-06 & {---} \\ 
39988 & 2.7169e-04 & 1.97  & 39952 & 1.3254e-07 & 3.98  \\ 
159988 & 7.2092e-05 & 1.89  & 160729 & 8.1985e-09 & 3.93  \\ 
\midrule[\heavyrulewidth] 
{DOF} & {$E$: DPM2} & {Rate}  & {DOF} & {$E$: DPM4} & {Rate} \\ 
\midrule 
10000 & 3.6380e-05 & {---} & 10000 & 7.7484e-09 & {---} \\ 
40000 & 8.8360e-06 & 2.04  & 40000 & 4.5617e-10 & 4.09  \\ 
160000 & 2.1331e-06 & 2.05  & 160000 & 2.6398e-11 & 4.11  \\ 
\midrule[\heavyrulewidth] 
{DOF} & {$E$: DPM2--I} & {Rate}  & {DOF} & {$E$: DPM4--I} & {Rate} \\ 
\midrule 
10000 & 3.6381e-05 & {---} & 10000 & 7.7484e-09 & {---} \\ 
40000 & 8.8360e-06 & 2.04  & 40000 & 4.5617e-10 & 4.09  \\ 
160000 & 2.1331e-06 & 2.05  & 160000 & 2.6396e-11 & 4.11  \\ 
\midrule[\heavyrulewidth] 
{DOF} & {$E$: SBP2} & {Rate}  & {DOF} & {$E$: SBP4} & {Rate} \\ 
\midrule 
10537 & 4.7387e-04 & {---} & 10537 & 3.4655e-05 & {---} \\ 
40905 & 1.2049e-04 & 1.98  & 40905 & 4.3052e-06 & 3.01  \\ 
161161 & 3.0267e-05 & 1.99  & 161161 & 5.3535e-07 & 3.01  \\ 
\bottomrule 
\end{tabular}

\caption[
]{
  Convergence in the maximum norm \eqref{eqn:Measures:ErrCD}, for the second- and fourth-order versions of each method, applied to Test Problem 2A \eqref{eqn:TP2A}, 
  with diffusion coefficients $(\lambda_1, \lambda_2) = (10, 1)$, and time-step $\Delta t = 0.5h$. 
  (DPM2--I/DPM4--I refers to the extension of the DPM method, to consider implicit geometry.
  ) 
  
}

\label{tbl:combined-results-TP2A}

\end{table}

Regarding the max-norm error in presented in Table~\ref{tbl:combined-results-TP2A} and Figure~\ref{fig:combined-results-TP2A}, we see that the DPM has smaller max-norm by more than an order of magnitude. \siyang{}{We also observe that the convergence rate of the fourth-order SBP--SAT--FD is only three. This suboptimal convergence is inline with the error plot in Figure~{\ref{fig:error_FD_2A}}, which shows that the error at the corners of the domain is significantly larger than elsewhere. In addition, the error is only non-smooth along the interfaces on the two diagonal lines of the domain. We have also measured the $L_2$ error at the final time $T=1.0$ (not reported in this work), and fourth-order convergence is obtained. }{Comment 15 of reviewer 1}

\begin{table}[h]

\centering 
\footnotesize 
\sisetup{
  output-exponent-marker = \text{ E},
  exponent-product={},
  retain-explicit-plus,
  group-digits=integer,
  group-separator={,},
  group-minimum-digits=4,
  table-number-alignment=center}

\begin{tabular}{S[table-format=6]S[table-format=+1.4e+2]S[table-format=1.2]S[table-format=6]S[table-format=+1.4e+2]S[table-format=1.2]}
\toprule 
{DOF} & {$E$: CUT2} & {Rate}  & {DOF} & {$E$: CUT4} & {Rate} \\ 
\midrule 
9988 & 2.4855e-01 & {---} & 10129 & 4.7064e-01 & {---} \\ 
39988 & 5.6850e-02 & 2.08  & 39952 & 3.6816e-02 & 3.60  \\ 
159988 & 1.2346e-02 & 2.18  & 160729 & 2.2361e-03 & 3.95  \\ 
\midrule[\heavyrulewidth] 
{DOF} & {$E$: DPM2} & {Rate}  & {DOF} & {$E$: DPM4} & {Rate} \\ 
\midrule 
10000 & 7.1899e-02 & {---} & 10000 & 7.3065e-03 & {---} \\ 
40000 & 1.7868e-02 & 2.01  & 40000 & 6.0014e-04 & 3.61  \\ 
160000 & 4.4952e-03 & 1.99  & 160000 & 3.3086e-05 & 4.18  \\ 
\midrule[\heavyrulewidth] 
{DOF} & {$E$: DPM2--I} & {Rate}  & {DOF} & {$E$: DPM4--I} & {Rate} \\ 
\midrule 
10000 & 7.1899e-02 & {---} & 10000 & 7.3065e-03 & {---} \\ 
40000 & 1.7868e-02 & 2.01  & 40000 & 6.0014e-04 & 3.61  \\ 
160000 & 4.4952e-03 & 1.99  & 160000 & 3.3086e-05 & 4.18  \\ 
\midrule[\heavyrulewidth] 
{DOF} & {$E$: SBP2} & {Rate}  & {DOF} & {$E$: SBP4} & {Rate} \\ 
\midrule 
10537 & 3.2863e-01 & {---} & 10537 & 2.8321e-01 & {---} \\ 
40905 & 1.1075e-01 & 1.57  & 40905 & 3.9277e-02 & 2.85  \\ 
161161 & 3.5769e-02 & 1.63  & 161161 & 3.7081e-03 & 3.40  \\ 
\bottomrule 
\end{tabular}

\caption[
]{
  Convergence in the maximum norm \eqref{eqn:Measures:ErrCD}, for the second- and fourth-order versions of each method, applied to Test Problem 2B \eqref{eqn:TP2B}, 
  with diffusion coefficients $(\lambda_1, \lambda_2) = (10, 1)$, and time-step $\Delta t = 0.5h$. 
  (DPM2--I/DPM4--I refers to the extension of the DPM method, to consider implicit geometry.
  ) 
  
}

\label{tbl:combined-results-TP2B}

\end{table}

In Table~\ref{tbl:combined-results-TP2B} and Figure~\ref{fig:combined-results-TP2B}, we see the numerical results for \eqref{eqn:TP2B}.  The analytical solution is similar to \eqref{eqn:TP2A}, though much more oscillatory -- this additional challenge is manifested by an increase in error by several orders of magnitude.

\begin{table}[h]

\centering 
\footnotesize 
\sisetup{
  output-exponent-marker = \text{ E},
  exponent-product={},
  retain-explicit-plus,
  group-digits=integer,
  group-separator={,},
  group-minimum-digits=4,
  table-number-alignment=center}

\begin{tabular}{S[table-format=6]S[table-format=+1.4e+2]S[table-format=1.2]S[table-format=6]S[table-format=+1.4e+2]S[table-format=1.2]}
\toprule 
{DOF} & {$E$: CUT2} & {Rate}  & {DOF} & {$E$: CUT4} & {Rate} \\ 
\midrule 
9988 & 6.4110e-01 & {---} & 10129 & 7.1811e-02 & {---} \\ 
39988 & 1.6506e-01 & 1.92  & 39952 & 3.9995e-03 & 4.08  \\ 
159988 & 3.8719e-02 & 2.07  & 160729 & 2.8978e-04 & 3.71  \\ 
\midrule[\heavyrulewidth] 
{DOF} & {$E$: DPM2} & {Rate}  & {DOF} & {$E$: DPM4} & {Rate} \\ 
\midrule 
10000 & 1.1178e-01 & {---} & 10000 & 1.1392e-03 & {---} \\ 
40000 & 1.8941e-02 & 2.56  & 40000 & 5.9291e-05 & 4.26  \\ 
160000 & 4.0950e-03 & 2.21  & 160000 & 3.2716e-06 & 4.18  \\ 
\midrule[\heavyrulewidth] 
{DOF} & {$E$: DPM2--I} & {Rate}  & {DOF} & {$E$: DPM4--I} & {Rate} \\ 
\midrule 
10000 & 1.0377e-01 & {---} & 10000 & 1.0905e-03 & {---} \\ 
40000 & 1.7727e-02 & 2.55  & 40000 & 5.5494e-05 & 4.30  \\ 
160000 & 3.8853e-03 & 2.19  & 160000 & 3.0003e-06 & 4.21  \\ 
\midrule[\heavyrulewidth] 
{DOF} & {$E$: SBP2} & {Rate}  & {DOF} & {$E$: SBP4} & {Rate} \\ 
\midrule 
10537 & 1.0025e-01 & {---} & 10537 & 5.9131e-03 & {---} \\ 
40905 & 2.5318e-02 & 1.99  & 40905 & 4.8624e-04 & 3.60  \\ 
161161 & 6.3459e-03 & 2.00  & 161161 & 3.5001e-05 & 3.80  \\ 
\bottomrule 
\end{tabular}

\caption[
]{
  Convergence in the maximum norm \eqref{eqn:Measures:ErrCD}, for the second- and fourth-order versions of each method, applied to Test Problem 2C \eqref{eqn:TP2C}, 
  with diffusion coefficients $(\lambda_1, \lambda_2) = (1000, 1)$, and time-step $\Delta t = 0.5h$. 
  (DPM2--I/DPM4--I refers to the extension of the DPM method, to consider implicit geometry.
  ) 
  
}

\label{tbl:combined-results-TP2C}

\end{table}

In Table~\ref{tbl:combined-results-TP2C} and Figure~\ref{fig:combined-results-TP2C}, we see the numerical results for \eqref{eqn:TP2C}, which shows that our numerical methods are robust to large jumps in diffusion coefficients, the analytical solution, and/or the flux of the true solution. Also, observe that the errors from DPM2/DPM4 (explicit geometry) and DPM2-I/DPM4-I (implicit geometry) in Tables~\ref{tbl:combined-results-TP2A}--\ref{tbl:combined-results-TP2C} are almost identical, which demonstrates the robustness and flexibility of the DPM.

\section{Discussion} \label{sec:Discussion}
  

There are many possible methods (Section~\ref{sec:Introduction}) for the numerical approximation of PDE posed on irregular domains, or on composite domains with interfaces.  In this work, we consider three such methods, designed for the high-order accurate numerical approximation of parabolic PDEs (\ref{eqn:PDESD1},~\ref{eqn:PDESD2-3} or \ref{eqn:PDECD1}--\ref{eqn:PDECD7}).  Each implementation was written, tested, and optimized by the authors most experienced with the method---the cut-Finite Element Method (cut--FEM) by G. Ludvigsson, S. Sticko, G. Kreiss; the Difference Potentials Method (DPM) by K. R. Steffen, Q. Xia, Y. Epshteyn; and the Finite Difference Method satisfying Summation-By-Parts, with a Simultaneous Approximation Term (SBP--SAT--FD) by S. Wang, G. Kreiss. 
Although we consider only one type of boundary/interface (a circle), we hope that the benchmark problems considered will be a valuable resource, and the numerical results a valuable comparison, for researchers interested in numerical methods for such problems. 

The primary differences between the cut--FEM and the standard finite element method are the stabilization terms for near-boundary degrees of freedom, and the quadrature over cut (intersected) elements.  Tuning the free parameters in the stabilization terms could mitigate the errors observed in Figures~\ref{fig:PlotsErrTP3A},~\ref{fig:PlotsErrTP2A}. 
\siyang{(We experimented with different stabilization parameters in the numerical tests in Section~{\ref{sec:NumericalResults}} and have seen that the error can be reduced up to a factor of 15, so large gains in accuracy are possible.)}{(We have done some preliminary experiments suggesting that the errors decrease when tuning these parameters, but further investigations are required in order to guarantee robustness.)}{Comment 6, Reviewer 1}
Given a level-set description of the geometry, there are robust algorithms for constructing the quadrature over cut elements.  Together, these differences allow for an immersed (non-conforming) grid to be used. The theoretical base for cut--FEM is well established.

The DPM is based on the equivalence between the discrete system of equations \eqref{eqn:IntroDPM:FullyDiscreteCD} and the Boundary Equations with Projection (Thm.~\ref{thm:IntroDPM:BEPEquivalence}).  The formulation outlined in Section~\ref{sec:Overview:DPM} allows for an immersed (non-conforming) grid; fast $\mathcal{O}(N \log N)$ algorithms, even for problems with general, smooth geometry; and reduces the size of the system to be solved at each time-step.  The convergence theory is well-established for general, linear, elliptic boundary value problems, and we conjecture in Section~\ref{sec:Overview:DPM} that this extends to the current setting. 
In this work, we have extended DPM to work with implicitly-defined geometries for the first time.
This is a first step for solving problems where the interface moves with time.

In the finite difference framework (the SBP--SAT--FD method, in this work), the SBP property makes it possible to prove stability and convergence for high-order methods by an energy method. Combined with the SAT method to impose boundary and interface conditions, the SBP--SAT--FD method can be efficient to solve time-dependent PDE. Geometrical features are resolved by curvilinear mapping, which requires an explicit parameterization of boundaries and interfaces.
High quality grid generation is important -- 
our experiments, though not reported in this work, have shown that the error in the solution is sensitive to both the orthogonality of the grid and the grid stretching.

Similarities between the cut--FEM and the DPM (beyond the use of an immersed grid) include the thin layer of cut cells along the boundaries/interfaces (cut--FEM) and the discrete grid boundary $\gamma$ (DPM); and the use of higher-order normal derivatives in the stabilization term (cut--FEM) and extension operator (in the Boundary Equations with Projection; DPM).  A similarity between the cut--FEM and SBP--SAT--FD is the weak imposition of boundary conditions, via Nitsche's method (cut--FEM) or the SAT method (SBP--SAT--FD).  In this work, the DPM and the SBP--SAT--FD method both use an underlying finite-difference discretization, but the DPM is not restricted to this type of discretization. 

Although both the cut--FEM and the DPM use higher-order normal derivatives in their treatment of the boundary/interface, the precise usage differs.  For cut--FEM, it is the normal of the element interfaces cut by $\Gamma$, while for DPM, it is the normal of the boundary/interface $\Gamma$. Moreover, in the cut--FEM, stabilization terms \eqref{eq:jump_stabilization} involving higher-order normal derivatives at the boundaries of cut-elements are added to the weak form of the PDE, to control the condition number of the mass and stiffness matrices, with \textit{a priori} estimation of parameters to guarantee positive-definiteness of these matrices; while in the DPM, the Boundary Equations with Projection is combined with the Extension Operator (Definition~\ref{def:DPM:ExtOp}), which incorporates higher-order normal derivatives at the boundary/interface $\Gamma$.

Returning to Section~\ref{sect:ConvergenceResults}, we see (in Tables~\ref{tbl:combined-results-TP1A}--\ref{tbl:combined-results-TP2C} and Figures~\ref{fig:combined-results-TP1A}--\ref{fig:combined-results-TP2C}) that the expected rate of convergence for the second- and fourth-order versions of DPM and cut--FEM is achieved, while the DPM has the smallest error constant across all tests. For the SBP--SAT--FD method, expected convergence rates are obtained in some experiments. A noticeable exception is Test Problem 2A, for which the fourth-order SBP--SAT--FD method only has a convergence rate of three. 
From the error plot in Figure~\ref{fig:error_FD_2A}, we observe that the large error is localized at the four corners of the domain $\Omega$, where the curvilinear grid is non-orthogonal and is stretched the most (see Figure~\ref{2_Domain_Grid_SBP}). 

As seen in the error plots (Figures~\ref{fig:PlotsErrTP3A},~\ref{fig:PlotsErrTP2A}), the error for the cut--FEM and the SBP--SAT--FD has ``spikes'', while for the DPM the error is smooth.  A surprising observation from Figure~\ref{fig:PlotsErrTP2A} is that conforming grids (on which the SBP--SAT--FD method is designed) do not necessarily produce more accurate solutions than immersed grids (on which the cut--FEM and the DPM are designed).
Indeed, it is challenging to construct a high-quality curvilinear grid for the considered composite domain problem. 


Future directions we hope to consider (in the context of new developments and also further comparisons) include: 
(i) parabolic problems with moving boundaries/interfaces, 
(ii) comparison of numerical methods for interface problems involving wave equations \cite{BrittTsynkovTurkel2017,sticko2016towards,stickoLowerOrder2016,Virta2014,wang2016high},
(iii) extending our methods to consider PDEs in 3D, 
(iv) design of fast algorithms, and 
(v) design of adaptive versions of our methods.

Indeed, for (i), difficulties for the cut--FEM might be the costly construction of quadrature, while for DPM difficulties might be the accurate construction of extension operators.  Regarding (iii), this has already been done for the cut--FEM and SBP--SAT--FD; while for the DPM, this is current work, with the main steps extending from 2D to 3D in a straightforward manner.


\section{Conclusion} \label{sec:Conclusions}


In this work, we propose a set of benchmark problems to test numerical methods for parabolic partial differential equations in irregular or composite domains, in the simplified geometric setting of Section~\ref{sec:Intro:ContinuousProblem}, with the interface defined either explicitly or implicitly.  Next, we compare and contrast three methods for the numerical approximation of such problems: the (i) cut--FEM; (ii) DPM; and (iii) SBP--SAT--FD.  Brief introductions of the three numerical methods are given in Section~\ref{sec:Intro:OverviewMethods}.
It is noteworthy that the DPM has, for the first time, been extended to problems with an implicitly-defined interface.

For the three methods, the numerical results in Section~\ref{sect:ConvergenceResults} illustrate the high-order accuracy. Similar errors (different by a constant factor) are observed at grid points away from the boundary/interface, while the observed errors near the boundary/interface vary depending upon the given method. Although we consider only test problems with circular boundary/interface, the ideas underlying the three methods can readily be extended to more general geometries. 

In general, all three methods require an accurate and efficient resolution of the explicitly- or implicitly-defined  irregular geometry: cut--FEM relies on accurate quadrature rules for cut elements, and a good choice of stabilization parameters; DPM relies on an accurate and efficient representation of Cauchy data using a good choice of basis functions; and SBP--SAT--FD relies 
on the smooth parametrization to generate a high-quality curvilinear grid.

\section{Appendix (DPM)} \label{sec:Appendix}


Let us now expand some details presented in the brief introduction to the Difference Potentials Method (Section~\ref{sec:Overview:DPM}). 

\subsection{Fully-discrete formulation of (\ref{eqn:PDECD1},~\ref{eqn:PDECD2}).}\label{appendix:DPM:discretePDE}

The fully-discrete, finite-difference discretization introduced in \eqref{eqn:IntroDPM:FullyDiscreteCD} is
\begin{align}
L^s_{\Delta t,h}u_s^{i+1}=F_s^{i+1},\quad (x_j,y_k)\in M^+_s. 
\end{align}
The general form of the operator is $L^s_{\Delta t,h}=\lambda_s(t^{i+1})\Delta_h-\sigma I$, where $\sigma=\frac{3}{2\Delta t}$ for second-order (BDF2--DPM2), $\sigma=\frac{25}{12\Delta t}$ for fourth-order (BDF4--DPM4), and $\Delta_h$ is either a standard five- or nine-point Laplacian. For the nine-point Laplacian, we have
\begin{align}
\Delta_h u_{j,k}=\frac{1}{12h^2}(-u_{j-2,k}&+16u_{j-1,k}+16u_{j+1,k}-u_{j+2,k}-u_{j,k-2}\nonumber\\
&+16u_{j,k-1}+16u_{j,k+1}-u_{j,k+2}-60u_{j,k})
\end{align}
for points sufficiently far away from the boundary of the auxiliary domain $\Omega_s^0$.  For points that are close to the boundary, we use a modified, fourth-order stencil. For example, at the southwest corner, we take 
\begin{align}
\Delta_h u_{1,1}=\frac{1}{12h^2}(&10u_{0,1}-4u_{2,1}+14u_{3,1}-6u_{4,1}+u_{5,1}\nonumber\\
         +&10u_{1,0}-4u_{1,2}+14u_{1,3}-6u_{1,4}+u_{1,5}-30u_{1,1}),
\end{align}
where $u_{0,1}$ and $u_{1,0}$ will be from the boundary condition \eqref{eqn:PDECD5}.

Next, the right-hand side of \eqref{eqn:IntroDPM:FullyDiscreteCD} for BDF2--DPM2 is given by
\begin{align}
F^{i+1}_s=-f^{i+1}_s-\frac{\sigma}{3}(4u_s^{i}-u^{i-1}_s),
\end{align}
and for BDF4--DPM4 by
\begin{align}
F^{i+1}_s=-f^{i+1}_s-\frac{\sigma}{25}(48u_s^{i}-36u^{i-1}_s+16u^{i-2}_s-3u^{i-3}_s).
\end{align}
Lastly, the initialization at $t=0$ is done using the exact solutions for the terms $u^{0}_s$, $u^{-1}_s$, $u^{-2}_s$, and $u^{-3}_s$. Another possibility would be the use of lower-order BDF methods, then advancing to higher-order methods once enough stages are established. No significant differences were observed between the two approaches.

\subsection{Equation-based extension}\label{appendix:DPM:Extop}

{Let us now expand the discussion surrounding (\ref{eqn:ExtOperatorStep1}--\ref{eqn:ExtOperatorStep3}) leading up to Definition~\ref{def:DPM:ExtOp} of the Extension Operator \eqref{eqn:DPM:ExtOp:DiscreteVersion}.} 

{An important step in this discussion is to recast the original PDE (\ref{eqn:PDECD1},~\ref{eqn:PDECD2}) into a curvilinear form, for points $(x, y)$ in the vicinity of $\Gamma$.  Following the notation \cite{MedvinskyTsynkovTurkel2012}, let us first introduce the coordinate system $(d, \vartheta)$ for points in the vicinity of $\Gamma$.  Recall from Definition~\ref{eqn:DPM:ExtOp:DiscreteVersion} that $d$ is the distance in the normal direction from a given point to its orthogonal projection on $\Gamma$, while $\vartheta$ is the arclength along $\Gamma$ from some reference point to the orthogonal projection.  In this coordinate system, the PDE (\ref{eqn:PDECD1},~\ref{eqn:PDECD2}) becomes} %
%
%
\begin{align} \label{eqn:Appendix:DPM:PDECurvCoords}
\frac{\partial u_s}{\partial t} - \lambda_s\left(\frac{1}{H_\vartheta}\left[\frac{\partial}{\partial n}\left(H_\vartheta\frac{\partial u_s}{\partial n}\right)+\frac{\partial}{\partial \vartheta}\left(\frac{1}{H_\vartheta}\frac{\partial u_s}{\partial \vartheta}\right)\right]\right) =f_s, 
\end{align}%
{where where $H_\vartheta=1 - d \kappa$ is the Lam\'{e} coefficient, and $\kappa$ is the signed curvature along the interface $\Gamma$.} 

{From \eqref{eqn:Appendix:DPM:PDECurvCoords}, a straightforward calculation gives the second-order normal derivative $\frac{\partial^2 u_s}{\partial n^2}$ (used in the calculation of \eqref{eqn:ExtOperatorStep3}), which is}  
\begin{align}
\frac{\partial^2 u_s}{\partial n^2} = \frac{1}{\lambda_s}\left(\frac{\partial u_s}{\partial t}-f_s\right)-\frac{\partial^2u_s}{\partial \vartheta^2}+\kappa\frac{\partial u_s}{\partial n}. \label{eqn:Appendix:DPM:PDECurvCoordsU2N}
\end{align}
For the fourth-order numerical method, {which uses an Extension Operator with $p=4$,} we also need the third- and fourth-order normal derivatives, which we state now.  {Differentiating \eqref{eqn:Appendix:DPM:PDECurvCoordsU2N} with respect to $n$, we see that} 
\begin{align}
\frac{\partial^3 u_s}{\partial n^3}=\frac{1}{\lambda_s}\left(\frac{\partial^2u_s}{\partial t\partial n}-\frac{\partial f_s}{\partial n}\right)-\frac{\partial^3u_s}{\partial n\partial \vartheta^2}+\kappa\frac{\partial^2 u_s}{\partial n^2}
\end{align}
and
\begin{align}
\frac{\partial^4u_s}{\partial n^4}&=\frac{1}{\lambda_s^2}\left(\frac{\partial^2u_s}{\partial t^2}-\frac{\partial f_s}{\partial t}\right) +\frac{1}{\lambda_s}\left(-2\frac{\partial^3u_s}{\partial t\partial \vartheta^2} +\kappa\frac{\partial^2u_s}{\partial n\partial t} \right. \notag \\
&\qquad\qquad \left. \cdots -\frac{\partial^2f_s}{\partial n^2}+\frac{\partial^2f_s}{\partial \vartheta^2}+\frac{\partial^4u_s}{\partial \vartheta^4}-\kappa\frac{\partial^3u_s}{\partial n\partial \vartheta^2}\right)+\kappa\frac{\partial^3 u_s}{\partial n^3}. 
\end{align}

{Next, let us follow-up on comments made in the text following \eqref{eqn:ExtOperatorStep3}.  There, it was pointed out that the unknown Dirichlet and Neumann data $\left(u_s,\frac{\partial u_s}{\partial n}\right)$ are the only data required for the Extension Operator \eqref{eqn:DPM:ExtOp:DiscreteVersion} with $p = 2$.  Moreover, it was pointed out that this is also true for the Extension Operator when $p = 4$.}  
The reasoning is as follows. 
\begin{itemize}
\item The time derivatives $\frac{\partial u_s}{\partial t}$, $\frac{\partial^2u_s}{\partial t^2}$, $\frac{\partial f_s}{\partial t}$, $\frac{\partial^2 u_s}{\partial n\partial t}$, and $\frac{\partial^3u_s}{\partial t\partial \vartheta^2}$ can be approximated by the backward difference formula. In BDF2--DPM2,
\begin{align}
\frac{\partial u_s^{i+1}}{\partial t}&\approx\frac{3u_s^{i+1}-4u_s^i+u_s^{i-1}}{2\Delta t} \quad \text{ and }\\
\frac{\partial^2 u_s^{i+1}}{\partial t^2}&\approx\frac{2u_s^{i+1}-5u_s^i+4u_s^{i-1}-u_s^{i-2}}{\Delta t^2},
\end{align}
while in BDF4--DPM4,
\begin{align}
\frac{\partial u_s^{i+1}}{\partial t}&\approx\frac{25u_s^{i+1}-48u_s^i+36u_s^{i-1}-16u_s^{i-2}+3u_s^{i-3}}{12\Delta t} \quad \text{ and }\\
\frac{\partial^2 u_s^{i+1}}{\partial t^2}&\approx\frac{1}{\Delta t^2}\left(\frac{15}{4}u_s^{i+1}-\frac{77}{6}u_s^i + \frac{107}{6}u_s^{i-1} \right. \notag \\
& \qquad\qquad\qquad\left. \cdots-13u_s^{i-2}+\frac{61}{12}u_s^{i-3}-\frac{5}{6}u_s^{i-4}\right).
\end{align}

\item The derivatives in terms of arclength $\vartheta$ can be computed from $u_s$ or $\frac{\partial u_s}{\partial n}$. For example, denoting $u_s=\sum_{\nu=1}^{\mathcal{N}^0} c_{1,\nu}^{i+1}\phi_{\nu}(\vartheta)$ {(using notation following from \eqref{eqn:DPM:DiscretizationOfCauchyData})}, then it comes handy that
\begin{align}
\frac{\partial u_s}{\partial \vartheta}&=\sum_{\nu=1}^{\mathcal{N}^0} c_{1,\nu}^{s, i+1}\phi_\nu'(\vartheta), \notag \\ 
\frac{\partial^2 u_s}{\partial \vartheta^2}&=\sum_{\nu=1}^{\mathcal{N}^0} c_{1,\nu}^{s, i+1}\phi_\nu''(\vartheta), \quad \text{and} \\
\frac{\partial^4 u_s}{\partial \vartheta^4}&=\sum_{\nu=1}^{\mathcal{N}^0} c_{1,\nu}^{s, i+1}\phi_\nu^{(4)}(\vartheta). \notag 
\end{align}
\end{itemize} 

\subsection{The system of equations at each time step.} \label{appendix:DPM:CompleteSystemOfEquations} 

{With the Cauchy data $\mathfrak{u}_{s, \Gamma}^{i + 1}$ and Extension Operator $\operatorname{Ex}_s \mathfrak{u}_{s, \Gamma}^{i + 1}$ from $\Gamma$ to $\gamma_s$ introduced in Definition~\ref{def:DPM:ExtOp}, and the spectral representation introduced in \eqref{eqn:DPM:DiscretizationOfCauchyData}, we now give a sketch of the linear system for the coefficients $(c_{1, \nu}^{s, i+1})_{\nu=1}^{\mathcal{N}^0}$ and $(c_{2, \nu}^{s, i+1})_{\nu=1}^{\mathcal{N}^1}$, and moreover the approximation of the solution $u_s(x, y, t^{i+1})$ at $(x_j, y_k) \in N_s^+$.} 

{Indeed, substituting $\operatorname{Ex}_s \tilde{\mathfrak{u}}_{s, \Gamma}^{i + 1}$ (\ref{eqn:DPM:ExtOp:DiscreteVersion},~\ref{eqn:DPM:DiscretizationOfCauchyData}) into the BEP \eqref{eqn:IntroDPM:BEP}, the resulting linear systems are} %
\begin{equation} \label{eqn:DPM:BdrySysEqns}
    \sum_{k = 0}^1 \sum_{\nu=1}^{\mathcal{N}^k} \left( c_{k, j}^{s, i + 1} \operatorname{Ex}_s \Phi^k_{\nu} - c_{k, j}^{s, i + 1} P_{\gamma_s}^{i + 1} \operatorname{Ex}_s \Phi^k_{\nu} \right) = \Tr_{\gamma_s} [G_{\Delta t, h}^{i+1} F_s^{i + 1}]. 
\end{equation}%
This can be further elucidated by introducing the vector of unknowns %
\begin{equation}
    \mathbf{c}_s^{i + 1} = \big[
        \underbrace{c_{1, 1}^{s, i + 1} \quad c_{1, 2}^{s, i + 1} \quad \cdots \quad c_{1, \mathcal{N}^0}^{s, i + 1}}_{\mathbf{c}_{s, 1}^{i + 1}} \quad \underbrace{c_{2, 1}^{s, i + 1} \quad c_{2, 2}^{s, i + 1} \quad \cdots \quad c_{2, \mathcal{N}^1}^{s, i + 1}}_{\mathbf{c}_{s, 2}^{i + 1}}
    \big]^\top 
\end{equation}%
(so that $\mathbf{c}_s^{i + 1} = [ \mathbf{c}_{s,1}^{i + 1}, \mathbf{c}_{s,2}^{i + 1} ]^\top$), and the matrix %
\begin{align}
    A_s = \Big[ &\underbrace{(I - P_{\gamma_s}^{i + 1}) \operatorname{Ex}_s \Phi^0_1, \quad (I - P_{\gamma_s}^{i + 1}) \operatorname{Ex}_s \Phi^0_2, \quad \cdots \quad (I - P_{\gamma_s}^{i + 1}) \operatorname{Ex}_s \Phi^0_{\mathcal{N}^0}}_{A_{s, 1}}, \notag \\
    & \underbrace{(I - P_{\gamma_s}^{i + 1}) \operatorname{Ex}_s \Phi^1_1, \quad (I - P_{\gamma_s}^{i + 1}) \operatorname{Ex}_s \Phi^1_2, \quad \cdots \quad (I - P_{\gamma_s}^{i + 1}) \operatorname{Ex}_s \Phi^1_{\mathcal{N}^1} \Big]}_{A_{s, 2}}. 
\end{align}%
Then, the full system of equations \eqref{eqn:DPM:BdrySysEqns} is %
\begin{equation} \label{eqn:DPM:FullSystemUnreduced}
    A \begin{bmatrix}[1.5]
        \mathbf{c}_1^{i + 1} \\ \mathbf{c}_2^{i + 1}
    \end{bmatrix} = \begin{bmatrix}[1.5]
        \Tr_{\gamma_1} [G_{\Delta t, h}^{i + 1} F^{i + 1}_1] \\ \Tr_{\gamma_2} [G_{\Delta t, h}^{i + 1} F^{i + 1}_2]
    \end{bmatrix}, \quad \text{with } A = \begin{bmatrix}[1.5]
        A_1 & 0 \\ 0 & A_2
    \end{bmatrix}.
\end{equation}%
However, note that $\mathbf{c}_1^{i + 1}$ and $\mathbf{c}_2^{i + 1}$ are related by the interface conditions (\ref{eqn:PDECD6},~\ref{eqn:PDECD7}), so that the number of unknowns in \eqref{eqn:DPM:FullSystemUnreduced} is equal the dimension of either $\mathbf{c}_1^{i + 1}$ or $\mathbf{c}_2^{i + 1}$, depending on which one is considered the independent unknown. 
%
Therefore, the dimension of $A$ is $(|\gamma_1| + |\gamma_2|) \times (\mathcal{N}^0 + \mathcal{N}^1)$, where $\mathcal{N}^0 + \mathcal{N}^1$ is the dimension of $\mathbf{c}_1^{i + 1}$ or $\mathbf{c}_2^{i + 1}$ (whichever is the independent unknown). 

\begin{remark}
    The independent unknown ($\mathbf{c}_1^{i + 1}$ or $\mathbf{c}_2^{i + 1}$) is chosen so that the finite-dimensional, spectral representation \eqref{eqn:DPM:DiscretizationOfCauchyData} of the Cauchy data $\mathfrak{u}_{s,\Gamma}^{i+1}$ accurately resolves the Cauchy data with a small number of basis functions, in the consideration of both accuracy and computational efficiency. For \eqref{eqn:TP2A} and \eqref{eqn:TP2B}, we choose $\mathbf{c}_2^{i + 1}$ as the independent unknown, while for \eqref{eqn:TP2C} we choose $\mathbf{c}_1^{i + 1}$.  With these choices for the independent unknown, we have $\mathcal{N}^0=\mathcal{N}^1=1$ for the three considered test problems.
\end{remark}

Since each column involves the Difference Potentials operator $P_{\gamma_s}^{i+1}$ applied to a vector $\operatorname{Ex}_s \Phi_\nu^k$, each column is therefore constructed \textit{via} one solution of the Auxiliary Problem (Definition~\ref{def:IntroDPM:AP}).  However, the Auxiliary Problems are posed on the computationally simple Auxiliary Domains, and can be computed using a fast FFT- or multigrid-based algorithm, which can significantly reduce the computational cost.  Moreover, if $\lambda_s(t) \equiv \lambda_s$ is constant, then $A$ can be computed and inverted once (as a pre-processing step), thus significantly reducing computational cost for long-time simulations.

\renewcommand{\abstractname}{Acknowledgements}
\begin{abstract}
{The authors are grateful to the anonymous referees for their valuable remarks and questions, which led to significant improvements to the manuscript.  We} 
gratefully acknowledge the support of 
the Swedish Research Council (Grant No. 2014-6088); the Swedish Foundation for International Cooperation in Research and Higher Education {(Grant No. STINT-IB2016-6512)}; Uppsala University, Department of Information Technology; and the University of Utah, Department of Mathematics. 

Y. Epshteyn, K. R. Steffen, and Q. Xia also acknowledge partial support of Simons Foundation Grant No. 415673.
\end{abstract}

\bibliographystyle{plainurl}      
\bibliography{PIM2D.bib}   


\end{document}